\documentclass[leqno,12pt]{article}
\usepackage{latexsym}
\usepackage{amssymb}
\usepackage{amsmath}
\usepackage{bbm}
\usepackage{color}
\usepackage{graphicx}
\usepackage{epsfig}
\usepackage[latin1]{inputenc}
\newtheorem{thm}{Theorem}[section]
\newtheorem{lem}[thm]{Lemma}
\newtheorem{prop}[thm]{Proposition}
\newtheorem{cor}[thm]{Corollary}

\newtheorem{rem}[thm]{Remark}
\newtheorem{conj}[thm]{Conjecture}

\DeclareMathOperator{\cotan}{cotan}
\DeclareMathOperator{\SC}{SC}
\DeclareMathOperator{\GW}{GW}
\DeclareMathOperator{\HP}{HP}
 \newcommand{\ii}{{\mathrm{i}}}
\newcommand{\lal}{\langle\langle}
\newcommand{\rar}{\rangle\rangle}
\newcommand{\V}{\mathcal V}
\newcommand{\abold}{\boldsymbol{\alpha}}
\newcommand{\rbold}{\boldsymbol{\rho}}
\newcommand{\gbold}{\boldsymbol{\gamma}}
\newcommand{\fbold}{\boldsymbol{\varphi}}
\newcommand{\Fbold}{\boldsymbol{\Phi}}
\newcommand{\bbold}{\boldsymbol{b}}
\newcommand{\cbold}{\boldsymbol{c}}

\makeatletter\@addtoreset{equation}{section}\makeatother
\title{Sum rules and large deviations for spectral measures\\ on the unit circle}
\author{{\small Fabrice Gamboa}\footnote{ Universit\'e Paul Sabatier, Institut de Math\'ematiques de Toulouse,  31062 Toulouse Cedex 9, France, 
gamboa@math.univ-toulouse.fr}
\and{\small Jan Nagel}\footnote{Technische Universitat M\"unchen, Fakult\"at f\"ur Mathematik, Boltzmannstr. 3, 85748 Garching, Germany,  e-mail: jan.nagel@tum.de}
\and{\small Alain Rouault}\footnote
{Laboratoire de Mathématiques de Versailles, UVSQ, CNRS, Université Paris-Saclay, 78035-Versailles Cedex France, e-mail: alain.rouault@uvsq.fr}}
\begin{document}
\maketitle


\newcommand {\eref}[1]{(\ref{#1})}
\newcommand{\ea}{\end{array}}
\newcommand{\beqohne}{\begin{eqnarray*}}
\newcommand{\eeqohne}{\end{eqnarray*}}
\newcommand{\beohne}{\begin{equation*}}
\newcommand{\eeohne}{\end{equation*}}
\newcommand{\R}{\mathbb{R}}
\newcommand{\N}{\mathbb{N}}
\newcommand{\T}{\mathbb{T}}
\newcommand{\E}{\mathbb{E}}
\newcommand{\1}{\mathbf 1}
\def\proof{\noindent{\bf Proof:}\hskip10pt}
\def\QED{\hfill\vrule height 1.5ex width 1.4ex depth -.1ex \vskip20pt}
\def \sur#1#2{\mathrel{\mathop{\kern 0pt#1}\limits^{#2}}}
\def \el{\sur{=}{(d)}}
\newcommand{\indi}{\mathbbm{1}} 

\newcommand{\la}{\lambda}
\newcommand{\lai}{\la_i}
\newcommand{\laj}{\la_j}
\newcommand{\mun}{\mu^{(n)}}
\newcommand{\muun}{\mu^{(n)}_{\u}}
\newcommand{\mutn}{\tilde{\mu}^{(n)}}

\newcommand{\ga}{\gamma}
\newcommand{\gai}{\ga_i}
\newcommand{\gaj}{\ga_j}
\def \u{{\tt u}}
\def \d{{\tt d}}
\def \g{{\tt g}}
\def \w{{\tt w}}
\def \m{{\tt m}}
\def \p{{\tt p}}
\def \be{\begin{eqnarray*}}
\def \ee{\end{eqnarray*}}
\def \ben{\begin{eqnarray}}
\def \een{\end{eqnarray}}
\newcommand{\tr}{\mathrm{tr}}
\newcommand{\Dir}{\mathrm{Dir}}
\newcommand{\Q}{\mathbb{Q}}
\newcommand{\Pnv}{\mathbb P^{(n)}_V}
\newcommand{\munI}{\mu^{(n)}_I}
\newcommand{\muunI}{\mu^{(n)}_{\u,I}}
\newcommand{\mutnI}{\tilde{\mu}^{(n)}_I}
\newcommand{\munIj}{\mu^{(n)}_{I(j)}}
\newcommand{\muunIj}{\mu^{(n)}_{\u,{I(j)}}}
\newcommand{\bmuunIj}{\bar{\mu}^{(n)}_{\u,{I(j)}}}
\newcommand{\mutnIj}{\tilde{\mu}^{(n)}_{I(j)}}
\newcommand{\sn}{^{(n)}}
\newcommand{\Ir}{\mathcal{J}}
\newcommand{\Fr}{\mathcal{F}}
\newcommand{\Sr}{\mathcal{S}}
\newcommand{\ap}{\alpha^+}
\newcommand{\am}{\alpha^-}
\newcommand{\bp}{b^+}
\newcommand{\bm}{b^-}
\newcommand{\si}{\sigma}
\newcommand{\siIj}{\sigma_{I_j}}

\newcommand{\lap}{\la^+}
\newcommand{\lam}{{\la^-}}
\newcommand{\lapm}{{\la^\pm}}
\newcommand{\lapj}{\la^+(j)}
\newcommand{\lamj}{\la^-(j)}
\newcommand{\lapjM}{\la_M^+(j)}
\newcommand{\lamjM}{\la_M^-(j)}
\newcommand{\lapmj}{\la^\pm(j)}
\newcommand{\tlap}{\tilde{\la}^+}
\newcommand{\tlam}{{\tilde{\la}^-}}
\newcommand{\tlapm}{{\tilde{\la}^\pm}}
\newcommand{\tlapj}{\tilde{\la}^+(j)}
\newcommand{\tlamj}{\tilde{\la}^-(j)}
\newcommand{\tlapmj}{\tilde{\la}^\pm(j)}

\newcommand{\gap}{\ga^+}
\newcommand{\gam}{\ga^-}
\newcommand{\bgap}{\bar{\ga}^+}
\newcommand{\bgam}{\bar{\ga}^-}
\newcommand{\gapm}{\ga^\pm}
\newcommand{\gapj}{\ga^+(j)}
\newcommand{\gamj}{\ga^-(j)}
\newcommand{\bgapj}{\bar{\ga}^+(j)}
\newcommand{\bgamj}{\bar{\ga}^-(j)}
\newcommand{\gapmj}{\ga^\pm(j)}
\newcommand{\tgap}{\tilde{\ga}^+}
\newcommand{\tgam}{\tilde{\ga}^-}
\newcommand{\tgapm}{\tilde{\ga}^\pm}
\newcommand{\tgapj}{\tilde{\ga}^+(j)}
\newcommand{\tgamj}{\tilde{\ga}^-(j)}
\newcommand{\tgapmj}{\tilde{\ga}^\pm(j)}

\newcommand{\Ep}{E^+}
\newcommand{\Em}{E^-}

\begin{abstract}
This work is a companion paper of \cite{GaNaRo} and \cite{GaNaRomat} (see also \cite{BSZ}). We continue to explore the connections between large deviations for random objects issued from random matrix theory and sum rules. Here, we are concerned essentially with measures on the unit circle whose support is an arc that is possibly proper. We particularly focus on two  matrix models. The first one is the Gross-Witten  ensemble. 
In the gapped regime we give a probabilistic interpretation of a Simon sum rule.  The second matrix model is the Hua-Pickrell ensemble.  Unlike the Gross-Witten ensemble the potential is here infinite at one point. Surprisingly, but as in \cite{GaNaRo}, we obtain a completely new sum rule for the deviation to the equilibrium measure of the Hua-Pickrell ensemble. The case of spectral matrix measures is also studied. Indeed, in the case of Hua-Pickrell ensemble, we extend our earlier works on large deviation  for spectral matrix measure \cite{GaNaRomat} and get here also a completely new sum rule.
\end{abstract}
\bigskip   

\smallskip

\maketitle

\section{Introduction}
Two of the most famous sum rules are Szeg\H{o}'s formula and the Killip-Simon sum rule. They are related to the theory of orthogonal polynomials on the unit circle (OPUC) and on the real line (OPRL),  respectively.

In the  OPUC frame,  the Szeg\H{o}-Verblunsky theorem (see \cite{Simon-newbook}, Theorem 1.8.6) concerns a  
deep relationship between 
the entropy of a measure $\mu$ supported by the unit circle
\begin{eqnarray*}
\mathbb T = \{ z \in \mathbb C : |z| = 1\} \simeq \{e^{\ii \theta} : \theta  \in [0, 2\pi)\}\,.
\end{eqnarray*}
and the coefficients involved in the construction of the orthogonal polynomial sequence in $L^2(\mu)$.
  More precisely, the 
recurrence relation  between two successive monic orthogonal polynomials $\phi_{k+1}$ and $\phi_{k}$ (where $\deg\phi_k=k$, $k\geq 0$) associated with a probability measure $\mu$ on the unit circle $\T$ supported by at least $k+1$ points 
involves a complex number $\alpha_k$
and may be written as
\begin{equation}
\label{recpolycirc}
\phi_{k+1}(z)=z\phi_{k}(z)-\overline{\alpha}_k\phi_{k}^*(z),\quad\mbox{ where } \quad \phi_{k}^*(z):=z^k\overline{\phi_k(1/\bar{z})}.
\end{equation}
The complex number $\alpha_k=-\overline{\phi_{k+1}(0)}$ is the so-called Verblunsky coefficient. In other contexts, it is also called Schur, Levinson, Szeg\H{o} coefficient or even canonical moment (\cite{DeSt97}). 
Let \[\mathbb D := \{ z \in \mathbb C : |z| < 1\}\] be the open unit disk. There are two different situations: 
when $\mu$ has a finite support of $n$ points, the coefficients satisfy $\alpha_k \in \mathbb D$ for $0 \leq k \leq n-2$ and $\alpha_{n-1} \in \mathbb T$, when $\mu$ has an infinite support, all the $\alpha_k$'s lie in $\mathbb D$. 

The Szeg\H{o}-Verblunsky  theorem is the identity
\begin{equation}
\label{segverth0}
\frac{1}{2\pi}\int_{0}^{2\pi}\log g_{\mu}(\theta)d\theta = \sum_{k\geq 0}\log(1-|\alpha_k|^2)\,,
\end{equation}
where the Lebesgue decomposition 
of $\mu$  
 is
$$d\mu(\theta)= g_{\mu}(\theta)\frac{d\theta}{2\pi}+d\mu_s(\theta)\,,$$
and where both sides of (\ref{segverth0}) are simultaneously finite or infinite. Changing the signs in both sides of this equation leads to 
\begin{equation}
\label{SVsum}
\mathcal K(\operatorname{UNIF} | \mu) = -\sum_{k \geq 0} \log(1 - |\alpha_k|^2)
\end{equation}
where, for probability measures $\nu$ and $\mu$, $\mathcal K (\nu | \mu)$ denotes the Kullback-Leibler divergence or relative entropy of $\nu$ with respect to $\mu$ (see (\ref{kukul})),  and $\operatorname{UNIF}$ is the normalized Lebesgue measure on  $\mathbb T$.

In the OPRL frame, for a probability measure $\mu$ having an infinite support, a.k.a. nontrivial case  (resp. with a finite support consisting of $n>0$ points, a.k.a. trivial case), the orthonormal polynomials associated to $\mu$  (with positive leading coefficients) obtained by applying the orthonormalizing Gram-Schmidt procedure 
to the sequence $1, x, x^2, \dots$ obey the 
recurrence relation
\begin{align} \label{polrecursion}
xp_k(x) = a_{k+1} p_{k+1}(x) + b_{k+1} p_k (x) + a_{k} p_{k-1}(x) 
\end{align}
for $ k \geq 0$ (resp. for $0 \leq k \leq n-1$). 
The Jacobi parameters $(a_k),(b_k)$ satisfy $b_k \in \mathbb R, a_k > 0$. Notice that here the orthogonal polynomials are not monic but normalized in $L^2(\mu)$.

To describe the Killip-Simon sum rule, we need some more  notations. Let $\mathcal{M}_1(I)$ denote the set of all probability measures on  $I$, a subset of $\mathbb{R}$ or of 
 $\mathbb{T}$. For $\alpha^-<\alpha^+$, let $\Sr_1^\mathbb{R}(\am,\ap)$ be the set of all probability measures $\mu$ on $\R$ with 
\begin{itemize}
\item[(i)] $\operatorname{supp}(\mu) = J \cup \{\lambda_i^-\}_{i=1}^{N^-} \cup \{\lambda_i^+\}_{i=1}^{N^+}$, where $J\subset [\am,\ap]$, $N^-,N^+\in\N_0\cup\{\infty\}$ and 
\begin{align*}
\lambda_1^-<\lambda_2^-<\dots <\am \quad \text{and} \quad \lambda_1^+>\lambda_2^+>\dots >\ap .
\end{align*}
\item[(ii)] If $N^-$ (resp. $N^+$) is infinite, then $\lambda_i^-$ converges towards $\am$ (resp. $\lambda_i^+$ converges to $\ap$).
\end{itemize}
Such a measure $\mu$ will be written as
\begin{align}\label{muinS0}
\mu = \mu_{|I} +  \sum_{i=1}^{N^+} \gamma_i^+ \delta_{\lambda_i^+} + \sum_{i=1}^{N^-} \gamma_i^- \delta_{\lambda_i^-}\,.
\end{align} 

The reference probability measure is now the semicircle law 
\begin{equation}
\label{SC0}
\operatorname{SC}(dx) = \frac{1}{2\pi}\sqrt{4-x^2}\!\ \mathbbm{1}_{[-2, 2]}(x)\!\ dx\,.
\end{equation} 
Additionally, we set
\begin{align*}
\mathcal{F}_{\SC}^+(x) :=  \begin{cases}   \displaystyle \int_2^x \sqrt{t^2-4}\!\ dt = \tfrac{x}{2} \sqrt{x^2-4} - 2 \log \left( \tfrac{x+\sqrt{x^2-4}}{2}\right)& \mbox{ if }  x \geq 2,\\
      \infty & \mbox{ otherwise}
\end{cases}
\end{align*} 
and  $\mathcal{F}_H^-(x):=\mathcal{F}_H^+(-x)$ for $x\in\R$. 

For a probability measure $\mu\in \mathcal S_1^\mathbb{R}(-2,2)$  with recursion coefficients $(a_k),\, (b_k)$ as in \eqref{polrecursion},  
the Killip-Simon sum rule is the following equation 
(see \cite{Simon-newbook}, Theorem 3.5.5):
\begin{align}
\label{KSsum}
{\mathcal K}(\operatorname{SC}\! |\!\ \mu) +  \sum_{i=1}^{N^+} {\mathcal F}^+_{\SC}(\lambda_i^+)  +  \sum_{i=1}^{N^-} {\mathcal F}^-_{\SC}(\lambda_i^-) = \sum_{k= 1}^\infty \big( \tfrac{1}{2} b_k^2 + G(a_k^2)\big) ,
\end{align}
where $G(x) = x - 1 - \log x$, and where both sides may be infinite simultaneously.

The common feature of formulas (\ref{SVsum}) and (\ref{KSsum}) is that they state equalities between non-negative functionals. We can consider them as equalities of two discrepancies. On  the left side it is the reverse relative entropy with respect to some reference probability measure plus possibly a contribution of the outlying point masses. On the right side it is a sum vanishing only when the {\it coefficients} involved are those of the reference probability measure.
An important consequence of such an equality are equivalent conditions for the finiteness of both sides,  one formulated in terms of Jacobi coefficients and the other as a spectral condition. In the words of Simon \cite{gem}, these are the \emph{gems} of spectral theory.

In \cite{gamboacanonical} and \cite{GaNaRo}, we revisited these results from a probabilistic point of view and gave a new proof based on large deviations (as we we will explain below). We also refer to the work of Breuer et al. \cite{BSZ} which enlightens    non-probabilists about  \cite{GaNaRo}, \cite{gamboacanonical}. This allowed us in the OPRL case to discover new sum rules, corresponding to the Marchenko-Pastur and Kesten-McKay measures, respectively. The main interesting feature of  (\ref{KSsum}) is the role played by the outliers of the measure $\mu$, i.e. its  discrete masses located out of the support of the reference measure. 

Coming back to the OPUC case,  in the Szeg\H{o}-Verblunsky theorem (\ref{segverth0}) there is no outlier 
since the reference probability measure is supported by the full unit circle $\mathbb T$. Nevertheless, 
  there are some very interesting probability measures supported by a proper arc. 
In this paper,  we study 
 sum rules for families of reference probability measures that are possibly supported by a proper arc of the unit circle. In particular, we prove a 
new sum rule (see Theorem \ref{sumruleHP}) concerning the reference probability measure 
${\HP}_\d$ (see (\ref{limmeas}), $\d$ is a positive parameter) that is supported by the a proper arc depending on $\d$. 
Up to our knowledge, Theorem \ref{sumruleHP} is completely new.

Our method for finding and showing a sum rule relies on the large deviations properties for a sequence of random measures built on random matrices. Let us give in a nutshell the scheme of our probabilistic method. We interpret the measure $\mu$ as the realization of a (random) spectral measure of a pair $(M, e)$ where $M$ is a random normal operator (unitary or Hermitian) and $e$ a fixed vector in a Hilbert space $\mathcal H$. 

Let assume that $\dim\mathcal H = n \geq 1$. 
 Then, $\mu$ is a discrete probability measure which can be encoded as 
\begin{equation}
\label{meumeu}
\mu = \sum_{k=1}^n \w_k \delta_{\lambda_k}\,.
\end{equation}
A classical assumption is the invariance by any unitary conjugations of the law of $M$. Under this assumption,  the joint density of $(\lambda_1, \dots, \lambda_n)$ is  proportional to the square of the Vandermonde determinant multiplied by the exponential of some potential. Furthermore,  the distribution of the weights $(\w_1, \cdots, \w_n)$ is uniform on the simplex. This allows in the cases studied in \cite{gamboacanonical} and \cite{GaNaRo}, with convenient assumptions on the potential, to show that the random measure defined in (\ref{meumeu}) satisfies, as $n$ grows, a large deviation principle (LDP). The  speed of the LDP is $n$ and  the rate function is the left hand side of (\ref{SVsum}) or (\ref{KSsum}) or a similar expression. In the sum rules, the right hand side is obtained as the rate function seeing the random measure as encoded by its Verblunsky (OPUC) or Jacobi (OPRL) coefficients.  Since a rate function is unique, the equality of both sides follows straightforwardly.

Of course, there is a natural way to travel in both directions from $\mathbb T$ to $\mathbb R$. This is the so-called Cayley transform. So that,  the results obtained for random measures on $\mathbb R$ may be carried to random measures on $\mathbb T$. Nevertheless,  the confinement assumption made on the potential in \cite{GaNaRo} is not always true in  some 
 interesting cases on $\mathbb T$.  Two examples are particularly representative and more or less emblematic of studies in OPUC and in equilibrium measures on $\mathbb T$. The first one is the Gross-Witten (GW) ensemble  (gapped/ungapped regime), corresponding to a  potential continuous on $\mathbb T$. The second one is the Hua-Pickrell (HP) ensemble, corresponding to a potential infinite at one point. 
Both are distributions on the set (group) $\mathbb U(n)$ of unitary $n \times n$ matrices, (see \cite{Pritsker1}).
In the HP case, 
the potential, carried  on the real line, satisfies the confinement assumption.
It is then possible to use the results of \cite{GaNaRo} to state directly a LDP  
for the spectral measure. 
Moreover, since the {\it deformed} Verblunsky coefficients (see the Sections \ref{suveve} and \ref{sec:deformedVer} for the definitions of these coefficients) are independent with known distributions given in \cite{BNR}, the coding with these coefficients gives rise to a LDP and by uniqueness, we conclude with a new sum rule. This method is robust enough to be extended to the matrix case.

In the GW case, the 
 potential on the real line  satisfies only a weak growth assumption and we cannot use previous results. Nevertheless, we may work directly on $\mathbb T$, copying the scheme of
 proof of the real case, looking carefully at the differences.  We do not have exponential tightness for the extremal eigenvalues anymore, but since the potential is finite everywhere, we take benefit of the compactness of $\mathbb T$. It was the secret of Hardy's method  
 \cite{hardy2012note} to get the LDP for empirical spectral distribution under the weak growth assumption.
Besides and to be complete, we revisit the gapped case, for which the LDP is a direct consequence of \cite{gamboacanonical} and we give some probabilistic evidence for the celebrated sum rule due to  Simon \cite{simon05} (Theorem 2.8.1 therein). Notice that while we were revising this paper we have been aware of a recent work of Breuer {\it et al.} \cite{BSZ1} where very general sum rules, including the complete GW case, are shown using the large deviations approach.

For $p\in\mathbb{N}$, a normalized matrix measure  $\Sigma$ of size $p\times p$   on $\mathbb{T}$ is a matrix of signed complex measures, such that for any Borel set 
$A\subset \mathbb{T}$, $\Sigma(A)= (\Sigma_{i,j}(A))$ is Hermitian, non-negative definite and  such that $\Sigma(\mathbb{T})=\mathbf 1$. Here, $\mathbf 1$ denotes the identity matrix. As in the case of classical scalar measures, taking care of the non commutativity, it is possible to build associated right or left orthogonal matrix polynomial sequences (see Section \ref{MOPUC} and \cite{damanik2008analytic}).  These sequences satisfy recurrence relations as (\ref{recpolycirc}) involving matrix Verblunsky coefficients. Moreover, a Szeg\"o-Verblunsky identity (like (\ref{segverth0})) holds
(see \cite{delsarte} and \cite{derevy}). In the fields of probability and statistics the matrix measures and the corresponding Szeg\"o-Verblunsky identity have important applications in the spectral analysis and linear prediction of multivariate time series (see the survey
\cite{bingham2012multivariate}). As in the univariate frame we will give a completely new sum rule for matrix measure in terms of its {\it deformed} matrix Verblunsky coefficients. Indeed, Theorem \ref{sumruleHPmatrix} relates these coefficients with a matrix version of the reversed Kullback-Leibler divergence when the reference is the Hua-Pickrell matrix equilibrium measure.

The paper is organised as follows. In the next section we give some necessary notations and assumptions. 
In Section \ref{sexample} we describe the two main matrix models studied. Section \ref{sLDDPP} is devoted to our large deviation results for random spectral measures. The sum rules obtained from large deviation considerations are settled in Section \ref{sumrule}. At the end of this section, we present some connections with the past and present literature on sum rules and gems in the unit circle framework.  Extensions to matrix spectral measures are developed in Section \ref{smamat}. This section begins by some reminders on orthogonal matrix recursions and matrix  Verblunsky coefficients. All technical proofs are postponed to the last section.

\section{Notations, assumptions and tools}
\label{stoutou}
\subsection{Two encodings of  a probability measure on $\mathbb T$}
If $U$ is a unitary operator on a Hilbert space $\mathcal H$ and $e$ is  a cyclic vector for $U$, the spectral measure of the pair $(U ,e)$ is the unique 
probability measure $\mu$ on $\mathbb T$ such that
\begin{equation}
\langle e, U^k e\rangle = \int_\mathbb T z^k d\mu(z) \ \ (k \in \mathbb Z)\,.
\label{defspectralmeasure}
\end{equation}
Actually, $\mu$ is a unitary invariant for $(U,e)$.  If the dimension of $\mathcal H$ is $n$ and $e$ is cyclic for $U$,  let $\lambda_1= e^{\ii\theta_1}, \ldots, \lambda_n= e^{\ii\theta_n}$ be the  eigenvalues of $U$ and let $\psi_1, \ldots, \psi_n$ be a system of orthonormal eigenvectors. The spectral measure of the pair $(U,e)$ is then
\begin{align}\label{spectralmeasure}
\mun =  \sum_{k=1}^n \w_k\delta_{\lambda_k}\,,
\end{align}
with $\w_k= |\langle \psi_k, e\rangle|^2$ and $\delta_a$ is the Dirac measure at $a$. This measure is a weighted version of the empirical eigenvalue distribution
\begin{align}\label{empiricallaw}
\muun = \frac{1}{n} \sum_{k=1}^n \delta_{\lambda_k} \,.
\end{align}
Another invariant is the CMV (or 5-diagonal) reduction of $U$. 
Let us now describe shortly the CMV mapping between 5-diagonal matrices and spectral measures.

We consider $n \times n$ 
matrices  corresponding to  measures supported by 
$n$ points (trivial case) and semi-infinite matrices corresponding to measures with bounded infinite support (non-trivial case). In the basis $(\chi_k)_{k \geq 0}$ obtained by orthonormalizing $1, z , z^{-1}, z^2, z^{-2}, \dots$,
 the 
linear transformation $f(z) \rightarrow zf(z)$ 
 in $L^2(\mu)$ is represented by the matrix
\ben
\label{favardfini}
\mathcal C_\mu = \begin{pmatrix} \bar \alpha_0&\bar \alpha_1\rho_0 &\rho_1\rho_0&0&0&\dots\\
\rho_0& -\bar\alpha_1\alpha_0&-\rho_1\alpha_0 &0 &0 &\dots\\
0&\bar\alpha_2\rho_1&-\bar\alpha_2\alpha_1&\bar\alpha_3\rho_2&\rho_3\rho_2&\dots\\
0&\rho_2\rho_1& -\rho_2\alpha_1& -\bar\alpha_3\alpha_2&-\rho_3\alpha_2&\dots\\
0&0&0&\bar\alpha_4\rho_3&-\bar\alpha_4\alpha_3& \dots\\
\dots&\dots&\dots&\dots&\dots&\dots
\end{pmatrix}
\een
with 
\begin{equation}
|\alpha_k| < 1 \ \hbox {and} \  \rho_k = \sqrt{1- |\alpha_k|^2}
\end{equation}
for every $k \geq 0$ in the non-trivial  case  and  for $0 \leq k \leq n-1$ in the trivial case, with $|\alpha_{n-1}| = 1$ \cite{CMV}. If the  measure is  supported by $n$ points, then the last line is
\begin{eqnarray*}
\begin{cases}
0\quad  \dots\quad 0 \qquad 0 \quad\qquad \bar\alpha_{2r}\rho_{2r-1} \quad\qquad  - \bar\alpha_{2r}\alpha_{2r-1} \ \ &\hbox{if} \ n=2r+1,\\
0\quad \dots\quad  0\quad  \rho_{2r}\rho_{2r-1}\quad  - \rho_{2r}\alpha_{2r-1}\quad  -\bar\alpha_{2r+1}\alpha_{2r}\ \ &\hbox{if} \ n=2r+2 , \ r\geq 0\,.
\end{cases}
\end{eqnarray*}

 Actually, there
is a one-to-one correspondence between such a matrix, called finite CMV matrix  and a finitely supported measure. 
If $\mathcal C$ is  a such a matrix, 
we can take the first vector of the canonical basis as the cyclic vector $e$. Let  $\mu$  be the spectral measure associated to the pair $(\mathcal C, e_1)$, then $\mathcal C$ represents the multiplication by $z$ in the basis $(\chi_k)$ of orthonormal polynomials
 associated to $\mu$ and $\mathcal C=\mathcal C_\mu$.

More generally, if $\mu$ is a non-trivial probability measure on $\mathbb T$, 
we may apply the same Gram-Schmidt procedure and consider the associated semi-infinite CMV matrix $\mathcal C_\mu$.
Notice that  now we have $|\alpha_k| < 1$ for every $k$. The mapping $\mu \mapsto \mathcal C_\mu$ (called here the CMV mapping) 
 is a one to one correspondence between probability measures on $\mathbb T$ having  infinite support and this kind of CMV matrices. This result is 
sometimes called Verblunsky-Favard's theorem (see \cite{Simon3}, p. 432).

\subsection{The Cayley transform, random matrices and invariant models}
\label{sec:RMmodels}

We will switch several times between $\mathbb R$ and $\mathbb T$ and between distributions of unitary and Hermitian matrices. There is a natural connection  between these two sets and also between these  two sets of matrices. This transformation is the so-called Cayley transform or stereographical projection. We follow here partly \cite{Benko} in its presentation. Let $\bar{\mathbb R} = \mathbb R \cup \{\infty\}$ be the compactified real line, which is topologically isomorphic to $\mathbb T$.
Let $\tau$ be the 
Cayley transform defined by:
\begin{align}
\label{Cayleysmall}
\zeta \in \mathbb T \setminus\{1\}& \mapsto \tau(\zeta): = \ii \frac{1+\zeta}{1-\zeta}\\
\notag
\tau(1) &= \infty  \\ \notag \\
\label{Cayleysmalli}
x \in \mathbb R &\mapsto \tau^{-1} (x) = \frac{x-\ii}{x+\ii}\\
\notag
\tau^{-1}(\infty) &= 1\,.
\end{align}
 It is clear that $\tau^{-1}$ is a homeomorphism from $\mathbb R$ onto $\mathbb T \setminus \{1\}$. 
Let us notice the important relations
\begin{equation}\label{tau-}|\tau^{-1}(x) - \tau^{-1}(y)| = \frac{2|x-y|}{\sqrt{1+x^2}\sqrt{1+y^2}}\ , \ |1-\tau^{-1}(x)| = \frac{2}{\sqrt{1+x^2}}\,,\end{equation}

and with angular coordinates
\begin{eqnarray*}
\zeta=e^{\ii \theta} \Longleftrightarrow x = \tau(\zeta) = - \cot \theta/2
\end{eqnarray*}
and
\begin{equation}
\label{diff}
d\theta = \frac{2dx}{1+x^2}\,.
\end{equation}
At the level of measures, we will consider three spaces. First $\mathcal M_1 (\mathbb R)$ and $\mathcal M_1 (\mathbb T)$ are the spaces of probability measures on $\mathbb{R}$ and $\mathbb{T}$, respectively, equipped with the topology of the weak convergence. 
Finally we need to use the set $\mathcal M_{\leq 1}(\mathbb{R})$ of subprobabilities on $\mathbb R$, equipped with the topology of vague convergence.  
Let us define the  mapping  $\hat\tau: \mathcal M_1 (\mathbb T) \to \mathcal M_{\leq 1}(\mathbb{R})$ with $\hat\tau(\nu)$ defined by
\begin{equation}
\label{hattau}
\int _{\mathbb R}f(x) d\hat\tau(\nu)(x) = \int_{\mathbb T \setminus \{\bf 1\}} f(\tau(\zeta)) d\nu(\zeta)\,.
\end{equation}
for any $f \in \mathcal C_0 (\mathbb R)$, continuous and vanishing at infinity. 
The mapping $\hat\tau$ is continuous if we equip $\mathcal M_{\leq 1}$ with the topology of vague convergence. Notice that $\mathcal M_1 (\mathbb T)$ and $\mathcal M_{\leq 1}(\mathbb R)$ are compact sets. We endow all these sets with the corresponding Borel $\sigma$-algebra.
The image of the uniform distribution on $\mathbb T$ is  the Cauchy distribution on $\mathbb R$. 

Let $\mathbb U(n)$ be the set of unitary $n\times n$ matrices and let $\mathrm{I}_n$ the identity $n \times n$ matrix. The Cayley transform induces a transformation 
from $\mathbb U(n)  \setminus \{\mathrm{I}_n\}$ onto $\mathcal H_n$, the set of Hermitian $n \times n$ matrices by 
\begin{equation}\label{Cayleybig}
M= \tau(U) :=\ii \frac{\mathrm{I}_n+U}{\mathrm{I}_n-U} \Longleftrightarrow U = \tau^{-1}(M) = \frac{M-\ii \mathrm{I}_n}{M+\ii \mathrm{I}_n}\,,
\end{equation}
in the sense of functional calculus. 
We denote by $\mathbb P\sn$  the normalized Haar measure on $\mathbb U(n)$.
It is classical (Weyl integration formula, \cite{Blower} Thm. 2.6.5) that under $\mathbb P\sn$ the array of eigenvalues has a density with respect to the Lebesgue measure $d\zeta_1 \dots d\zeta_n$ on $\mathbb T^n$  which is proportional to
\[\left|\Delta(\zeta_1, \dots, \zeta_n)\right|^2 \, ,\]
where $\Delta$ is the Vandermonde determinant. 
More generally, it is usual to equip $\mathbb{U}(n)$ with a probability measure of the form
\begin{align}
\label{PnV}
d\mathbb P\sn_\V (U) = \frac{1}{\mathcal Z_n^\V} e^{-n \tr \V(U)} d\mathbb P\sn (U)\,,
\end{align}
where $\V$ satisfies a convenient integrability assumption and  $\mathcal Z_n^{\mathcal V}$ is  the normalizing constant.
The density of eigenvalues under $\mathbb P_\V\sn$ is then proportional to 
\begin{equation}
\label{VdMC}
\left|\Delta(\zeta_1, \cdots, \zeta_n)\right|^2 \exp \left( -n \sum_{i=1}^n \V(\zeta_i)\right)\, .\end{equation}
If $\mathbb Q\sn$ is the Haar measure on the additive group $\mathcal H_n$ of Hermitian matrices defined by
\[d \mathbb Q\sn (M) = \prod_{k=1}^n  dM_{kk} \prod_{1\leq k < l\leq n} d(\Re M_{k l}) \!\  \prod_{1\leq k < l\leq n} d(\Im M_{kl})\,,\]
the pushforward of $\mathbb P\sn$ by $\tau$ is the Cauchy ensemble whose density with respect to $\mathbb Q\sn$ is proportional to
$\det (\mathrm{I}_n + M^2)^{-n}$. 
Let us compute the density of the (real) eigenvalues of $M=\tau(U)$, which are the pushforward of the eigenvalues of $U$ by $\tau$, when $\mathbb U (n)$ is equipped with $\mathbb P^{(n)}_{\mathcal V}$.
From \eqref{tau-} we have, if $\zeta_i \not= 1$ for $ i \leq n$,
\[\left|\Delta(\zeta_1, \cdots, \zeta_n)\right| = 2^{n(n-1)/2} \left|\Delta(x_1, \cdots, x_n)\right| \prod_{i=1}^n (1+x_i^2)^{-(n-1)/2}\,,\]
and with (\ref{diff}) we conclude that the array of eigenvalues of $M$ has 
a joint density  proportional to
\[\left|\Delta(x_1, \cdots, x_n)\right|^2 \exp \left(- n\sum_{i=1}^n  V(x_i) \right) \, ,\]
with respect to the Lebesgue measure on $\mathbb R^n$, where the potentials $V$ and $\mathcal{V}$ are related by
\begin{equation}
\label{VtoV}
 V(x) = \V(\tau^{-1}(x)) + \log (1+x^2)\,.
\end{equation}
The inverse relation is 
\[\V(e^{\ii \theta}) = V(-\cotan \theta/2) +  \log |\sin \theta/2|\, .\] 
Of course, the same distribution of eigenvalues can be obtained by observing that the pushforward of \eqref{PnV} by $\tau$ is
\begin{align}
\label{QnV}
d\mathbb P\sn_V (M) = \frac{1}{\mathcal{Z}_n^V} e^{-n \tr \V(\tau^{-1}(M))} \det (\mathrm{I}_n + M^2)^{-n} d\mathbb Q\sn (M)\, .
\end{align}
Besides, it is known that in all these unitary invariant models, the matrix $[\psi_1, \dots, \psi_n]$ of eigenvectors (defined up to multiplication of each vector by a phase) is Haar distributed on $\mathbb U(n)$. In particular, the array of  weights $(\w_1,  \dots, \w_n)$ defined in \eqref{spectralmeasure} is uniformly distributed on the simplex $\{(\w_1, \dots, \w_n) \in [0,1]^n : \w_1+\dots + \w_n =1\}$.

Looking at the spectral measures, the above models can be generalized to log-gases. In this framework,   
 $n$ is the number of particles (or eigenvalues) denoted by $\lambda_1,\dots ,\lambda_n$, with the joint distributions $\Pi_\V\sn$ on $\mathbb T^n$ having  the 
 density
\begin{align}\label{generaldensity}
\frac{d\Pi_\V\sn(\lambda)}{d\lambda} = \frac{1}{Z_\V^n}
e^{- n\beta' \sum_{i=1}^n \V(\lambda_i)}\prod_{1\leq  i < j\leq n} |\lambda_i - \lambda_j|^\beta\,,
\end{align}
with respect to the Lebesgue measure $d\lambda = d\lambda_1 \cdots d\lambda_n$. Here $\beta' = \beta/2$ and $\beta>0$ is a parameter interpreted as the inverse temperature. 
Then it is possible to consider the CMV matrices having these particles as eigenvalues and weights distributed according to the density proportional to
\[\prod_{k=1}^n w_k^{\beta' -1}\]
with respect to the uniform measure on the simplex (the Dirichlet distribution of parameter $\beta'$). The correspondence \eqref{VtoV} between the potentials is now more complicated: the pushforward by $\tau$ gives the relation 
\begin{align}\label{VtoVbeta}
 V(x) = \V(\tau^{-1}(x)) + \left(1-\frac{1}{n}+\frac{1}{\beta' n} \right) \log (1+x^2)\, ,
\end{align}
that is, either $V$ or $\V$ is $n$-dependent. It is possible to treat this general case, see Remark 3.5 in \cite{GaNaRomat}, but for simplicity, we consider here only the case $\beta=2$.

\subsection{Assumptions on the potentials}

\subsubsection{Real line} 

We will assume that the potentials $V$ on $\mathbb R$ are finite and continuous everywhere. The classical assumption on the growth of the potential is

(R1s) Strong growth: 
\begin{equation}\label{strong}
\liminf_{|x|\to \infty} \frac{V(x)}{2\log|x|} > 1 \,.
\end{equation}
Recently, Hardy \cite{hardy2012note} introduced the weaker assumption

(R1w) Weak growth: 
\begin{equation}\label{HardyR}
\liminf_{|x|\rightarrow \infty}\!\  \left( V(x) - 2\log |x|\right) > - \infty\,.
\end{equation}
Under (R1w), the empirical distribution $\muun$ of eigenvalues $\lambda_1,\dots ,\lambda_n$ has a 
limit $\mu_V$ (in probability), which is the unique minimizer of 
\begin{align}
\label{ratemuu}
 \mu  \mapsto \mathcal E_V (\mu) := \int_{\mathbb R} V(x) d\mu(x) - \iint_{\mathbb R^2} \log |x-y| d\mu(x)d\mu(y)\,, \qquad \mu \in\mathcal{M}_1(\R).
\end{align}
The minimal value is denoted as
\[F_V = \mathcal E_V (\mu_V)\,.\]
Under (R1s), it is known that the support of $\mu_V$ is compact.
We will make in this case the additional assumption

\begin{itemize}
\item[(R2)] One-cut regime: the support of $\mu_V$ is a single compact interval $[\alpha^-, \alpha^+]$ 
( $\alpha^-< \alpha^+$).
\end{itemize}

The minimizer $\mu_V$  is characterized by the Euler-Lagrange variational equations
\begin{equation}
\label{ELR}
\Ir_V(x)\begin{cases} = 2\xi_V & \hbox{if}\ x \in [\alpha^-, \alpha^+]\\
\geq 2 \xi_V & \hbox{if} \ x \notin [\alpha^-, \alpha^+]
\end{cases}
\end{equation}
where $\Ir_V$ is the effective potential
\begin{align}
\label{poteffR}
\Ir_V (x) := V(x) -2\int_{\mathbb R} \log |x-\xi|\!\ d\mu_V(\xi)\,,
\end{align}
and $\xi_V$ is the so-called modified Robin constant. We will make use of the following assumption

\begin{itemize}
\item[(R3)] Control (of large deviations): 

$\Ir_V$ achieves its global minimum value on the complement of 
 $(\alpha^-, \alpha^+)$  only on the boundary of this set.
\end{itemize}

Furthermore, to obtain a non-variational expression for the rate we need the following conditions:
\begin{itemize}
\item[(R4)] Offcriticality: 
\begin{align*}
d\mu_V(x) = \frac{1}{\pi}S(x) \sqrt{(\alpha^+ -x) (x-\alpha^-)} 
\!\ dx
\end{align*}
where $S> 0$ on $[\alpha^-, \alpha^+]$. 
\item[(R5)] Analyticity: $V$ can be extended as a holomorphic function is some open neighborhood of $[\alpha^-, \alpha^+]$.
\end{itemize}
We remark that for $V$ strictly convex, the assumptions (R2), (R3) and (R4) are fulfilled (see \cite{Borot} and \cite{johansson1998fluctuations}). 
The following proposition is more or less classical, it follows for example from \cite{Deiftuniform} (proof of Theorem 3.6) or \cite{Albeverio} (Equation (1.13)).

\begin{prop}\label{effpotR}
If the conditions (R1s), and (R2) to (R5) are satisfied, then
\begin{equation}
\Ir_V(x) - 2 \xi_V = \begin{cases}\int_{\alpha^+}^x S(t) \sqrt{(t-\alpha^-)(t-\alpha^+)}\,\ dt & \hbox{if} \ x > \alpha^+\,,\\
\int_x^{\alpha^-} S(t) \sqrt{(\alpha^- -t)(\alpha^+ -t)}\,\ dt & \hbox{if} \ x < \alpha^-\,.
\end{cases}
\end{equation}
\end{prop}

\subsubsection{Unit circle}
Let $\varphi, \psi \in [0, 2\pi]$ be two angles with $\varphi<\psi$. We define $\widehat{[\varphi, \psi]}$ to be the arc $[e^{\ii \varphi}, e^{\ii \psi}] \subset \mathbb T$ where we go from $e^{\ii \varphi}$ to $e^{\ii \psi}$ in a counterclockwise direction. The potential $\V$ is supposed to be continuous on $\mathbb T \setminus\{\1\}$. We make the additional assumption: 
\begin{itemize}
\item[(T1)] $\V$ is  lower semicontinuous in $1$. Without loss of generality we may assume that \[\V(1) = \liminf_{z \to 1}  \V(z) \in (-\infty,\infty] \,.\]
\end{itemize}
This implies that there is a unique minimizer $\mu_\V$ of 
\begin{equation}
\label{4.70}
\mu \mapsto \mathcal E_\V (\mu) = \int_{\mathbb T} \V (z) d\mu(z) - \iint_{\mathbb T^2} \log|z-\zeta| \!\ d\mu(z)d\mu(\zeta)\,,\qquad \mu \in \mathcal{M}_1(\T).
\end{equation}
The minimal value is denoted by
\begin{equation}
F_\V = \mathcal E_\V(\mu_\V)\,.
\end{equation}
We will suppose that either the support of $\mu_\V$ is $\mathbb T$ or
\begin{itemize} 
\item[(T2)] One-cut regime: the support of $\mu_\V$ is a single arc $\widehat{[\alpha^-, \alpha^+]} \subset \widehat{(0, 2\pi)}$. 
\end{itemize}
In this case, $\mu_\V$ is characterized by the Euler-Lagrange equations:
\begin{equation}
\label{ELT}
\Ir_\V(e^{\ii \theta})
\begin{cases} = 2\xi_\V & \hbox{if}\ \theta \in  \widehat{[\alpha^-, \alpha^+]}\,,\\
\geq 2 \xi_\V & \hbox{if} \ \theta \notin  \widehat{[\alpha^-, \alpha^+]}\,, 
\end{cases}
\end{equation}
where $\Ir_\V$ is the effective potential
\begin{align}
\label{poteffT}
\Ir_\V (e^{\ii \theta}) := \V(e^{\ii \theta}) -2\int_{\mathbb T} \log |e^{\ii \theta}-\zeta|\!\ d\mu_\V(\zeta)\,,
\end{align}
and $\xi_\V$ is the modified Robin constant. Like in the case of the real line, we make the assumption
\begin{itemize}
\item[(T3)] Control (of large deviations): 

$\mathcal J_\V$ achieves its global minimum value on the complement of  $\widehat{[\alpha^-, \alpha^+]}$   only on the boundary of this set.
\end{itemize}

When $\theta \mapsto v(\theta) := \V(e^{\ii \theta})$ is convex, this condition is satisfied. Indeed, it is 
\begin{equation}
\label{cotg}
\int_{\mathbb T} \log |e^{\ii \theta}-\zeta|\!\ d\mu_\V(\zeta) = \int_{e^{\ii \alpha^-}}^{e^{\ii \alpha^+}} \log \left|\sin \frac{\theta - \varphi}{2}\right|\!\ d\mu_\V (e^{\ii \varphi}) + \log 2
\end{equation}
so that, for $0 < \theta < \alpha^-$, the function $\theta \mapsto \mathcal J_\V(e^{\ii \theta})$ is strictly convex, nonegative and vanishes for $\theta = \alpha^-$, hence is positive on $[0, \theta^-)$.
An analogous argument can be made (mutatis mutandis) for $\alpha^+ <\theta< 2\pi$.

The additional assumptions to obtain a non-variational expression for the rate are on the unit circle:
\begin{itemize}
\item[(T4)] Offcriticality: 
\begin{align*}
d\mu_\V(z) = \frac{1}{\pi}S(e^{\ii \theta}) \sqrt{|(e^{\ii \theta}-e^{\ii \alpha^-})(e^{\ii \theta}- e^{\ii \alpha^+})|}
\!\ d\theta
\end{align*}
where $S(e^{\ii \theta})> 0$ for  $\theta \in [\alpha^-, \alpha^+]$. 

\item[(T5)] Analyticity: $\V$ can be extended as a holomorphic function is some open neighbourhood in $\mathbb C$ of the arc  $\widehat{[\alpha^-, \alpha^+]}$.
\end{itemize}

\begin{rem} \label{rem:equivassump}
Assumption (T1) is equivalent via \eqref{Cayleysmall} and \eqref{VtoVbeta} to 
Hardy's assumption \eqref{HardyR}.
\end{rem}

Using the Cayley transform, we may carry the statement of Proposition \ref{effpotR} to the circle, taking into account that
\begin{align*}
\mathcal J_V(x)= \mathcal J_\V(e^{\ii \theta})+ 2 \int \log |1-\zeta | d\mu_\mathcal{V}(\zeta)
\end{align*}
with $x=\tau (e^{\ii \theta})$ and that
\begin{align*}
\sqrt{(x-\tau (e^{\ii \alpha^-}))(\tau (e^{\ii \alpha^+})-x))} = 2 \frac{\sqrt{|e^{\ii \theta}-e^{\ii \alpha^-}||e^{\ii \alpha^+}-e^{\ii \theta}|}}{|1 - e^{\ii\theta}|\sqrt{|1 - e^{\ii\alpha^-}||1-e^{\ii\alpha^+}|}} .
\end{align*}
This results in the following proposition.

\begin{prop}\label{effpotT}
If $\V$ satisfies assumptions (T1) to (T5), then
\begin{equation}
\mathcal J_\V(e^{\ii \theta}) - 2 \xi_\V = \begin{cases}
\int_\theta^{\alpha^-} S(e^{\ii \tau}) \sqrt{|(e^{\ii \tau}-e^{\ii \alpha^-})(e^{\ii \tau}- e^{\ii \alpha^+})|}\, d\tau & \hbox{if} \ \theta \in (0, \alpha^-]\,,\\
\int_{\alpha^+}^\theta S(e^{\ii \tau}) \sqrt{|(e^{\ii \tau}-e^{\ii \alpha^-})(e^{\ii \tau}- e^{\ii \alpha^+})|}\, d\tau & \hbox{if} \ \theta \in [\alpha^+, 2\pi)\,.
\end{cases}
\end{equation}
\end{prop}
\subsection{Large deviations}

\label{sular}

\subsubsection{Introduction}
In order to be self-contained, let us recall the definition of a large deviation principle. For a general reference of large deviation statements we refer to the book of \cite{demboz98} or to the Appendix D of \cite{agz}.

Let $U$ be a topological Hausdorff space with Borel $\sigma$-algebra $\mathcal{B}(U)$. We say that a sequence $(P_{n})$ of probability measures on $(U,\mathcal{B}(U))$ satisfies the large deviation principle (LDP) with speed $a_n$ and
rate function $\mathcal{I} : U \rightarrow [0, \infty]$  if:
\begin{itemize}
\item [(i)] $\mathcal I$ is lower semicontinuous.
\item[(ii)] For all closed sets $F \subset U$: $\displaystyle \qquad 
\limsup_{n\rightarrow\infty} \frac{1}{a_n} \log P_{n}(F)\leq -\inf_{x\in F}\mathcal{I}(x) $
\item[(iii)] For all open sets $O \subset U$: $\displaystyle \qquad 
\liminf_{n\rightarrow\infty} \frac{1}{a_n} \log P_{n}(O)\geq -\inf_{x\in O}\mathcal{I}(x) $
\end{itemize}
The rate function $\mathcal{I}$ is good if its level sets
$\{x\in U |\ \mathcal{I}(x)\leq a\}$ are compact for all $a\geq 0$. 
If in the conditions above, we replace {\it closed sets} by {\it compact sets}, we say that $(P_{n})$ satisfies a weak LDP. In this case, we can recover a LDP if the additional condition of exponential tighness is fulfilled:

For every $M > 0$ there exists a compact set $K_M \subset U$ such that
\[\limsup_{n\rightarrow\infty}  \frac{1}{a_n} \log P_{n}(U \setminus K_M) \leq -M\,.\]
In our case, the measures $P_n$ will be the distributions of the random spectral measures $\mu_n$ and we will say that the sequence of measures $\mu_n$ satisfies a LDP. All LDPs for spectral measures in this section are in the weak topology.

\subsubsection{LDP for ESD}

The most famous LDP in random matrix theory 
concerns
 the sequence of empirical spectral measures (ESD) as defined in \eqref{empiricallaw}. The improved version (in the case $\beta=2$) is
 
\begin{prop}[Hardy \cite{hardy2012note} Thm. 1.1]
\label{LDPmuuR}
If the potential $V$ in $\mathbb R$ satisfies assumption (R1w), 
then under $\Pnv$, the sequence of random probability measures $(\mu_\u\sn)$
satisfies in $\mathcal M_1(\mathbb R)$ the LDP with speed $n^2$ and good rate function
\begin{align*}I_V(\mu) := \mathcal E_V (\mu) - F_V
\end{align*} where $\mathcal E_V$  is defined  in (\ref{ratemuu})\,. \end{prop}

An equivalent statement may be claimed via the Cayley transform (see also Remark 2.4 in \cite{hardy2012note}).

\begin{cor}
\label{LDPmuuU}
If the potential $\V$ satisfies assumption (T1), 
then under $\mathbb P^{(n)}_{\mathcal V}$ (see \eqref{PnV}), the sequence of random probability measures $(\mu_\u\sn)$
satisfies in $\mathcal M_1(\mathbb T)$ the LDP with speed $n^2$ and good rate function
\begin{align*}I_\V( \mu) := \mathcal E_\V (\mu) - F_\V
\end{align*}
where $\mathcal E_\V$  is defined  in (\ref{4.70})\,. \end{cor}

\section{Our two main examples of matrix ensembles}
\label{sexample}

\subsection{Gross-Witten ensemble}
Let us consider the  Gross-Witten measure on $\mathbb U(n)$, absolutely continuous with respect to the Haar measure $\mathbb P\sn$, with density:
\begin{equation}
\label{GW0}
\frac{d\mathbb G\mathbb W_{\g}\sn}{d\mathbb P\sn}(U) := \frac{1}{{\mathcal Z}_n(\g)} \exp \left[\frac{n\g}{2}\tr (U + U^\dagger)\right] \,,
\end{equation}
where $\g\in \mathbb R$, $\mathcal Z_n (\g)$ is the normalizing constant and $U^\dagger$ is the adjoint of $U$.
For details and applications of this distribution we refer to \cite{HiaiP} p. 203, \cite{GrossW}, \cite{Wadia}. 
It is important   in the analysis of problems involving random permutations since (Gessel relation)
\[\mathcal Z_n\left(\frac{2\sqrt \lambda}{n}\right) = e^{\lambda} \mathbb P(\ell_{N_\lambda} =n)\]
where $N_\lambda$ is a Poisson random variable of parameter $\lambda$ and $\ell_N$ is the length of the longest increasing subsequence  of a random permutation of $\{1, \dots, N\}$ chosen uniformly (formula (1.14) in \cite{DJ}). 

The potential is 
\begin{equation}
\label{potGW}
\V_\g(z) =  - \g \Re (z) \,.
\end{equation}
Although the potential $\V_\g$ is not a convex function of $\theta$, it is known that for this example (T3) is satisfied, see Lemma 4.3 in \cite{DJ}.  
If $|\g| \leq 1$ (ungapped or strongly coupled phase), the equilibrium measure $\GW_\g$ is supported by $\mathbb T$:
\begin{equation}
\label{GW-}
\GW_\g(dz)
= \frac{1}{2\pi} (1 + \g \cos \theta)\!\ d\theta,\; (z = e^{\ii \theta } , \theta \in [-\pi, \pi)).
\end{equation}
Moreover, we have
\begin{eqnarray}
\label{cstGW-}
F^{GW}_\g &=&\g^2/2\,,\\
\label{xiGW-}
\xi^{GW}_\g &=&\g^2/4\,.
\end{eqnarray}

Let us recall from Simon \cite{simon05}, p. 86 that the equilibrium measure has Verblunsky coefficients
\begin{equation}
\label{alphalimGW}
\alpha_n(\GW_\g) = \begin{cases}\displaystyle  -\frac{x_+ - x_-}{x_+^{n+2} - x_-^{n+2}} & \hbox{if} \ |\g| < 1\\
\displaystyle \frac{(-\g)^{n+1}}{n+2}& \hbox{if} \ |\g| = 1\,,
\end{cases}
\end{equation}
%
where $x_\pm = -\g^{-1} \pm \sqrt{\g^{-2}-1}$
 are roots of the equation 
\[x + \frac{1}{x} = -\frac{2}{\g}\,.\]
We remark that the distribution $ \GW_\g$ has only nontrivial moments of order  $\pm 1$.

For $|\g| > 1$  (gapped or weakly coupled phase), let $\theta_{g}  \in [0, \pi]$ be such that \begin{equation}\label{eqthetag}\sin^2 (\tfrac{\theta_\g}{2}) = |\g|^{-1}\,.\end{equation}
For $\g > 1$, the equilibrium measure is 
\begin{equation}\label{GWeq}
\GW_\g (dz) = \frac{\g}{\pi}\cos(\tfrac{\theta}{2})\!\  \sqrt{\sin^2(\tfrac{\theta_\g}{2})-\sin^2(\tfrac{\theta}{2})}\!\ 1_{[-\theta_\g, \theta_\g]}\!\ d\theta\;, (z = e^{\ii \theta} , \theta \in [-\pi, \pi)).\end{equation}
Moreover, the free energy and the modified Robin constant are in the gapped case
\begin{eqnarray}
\label{cstGW}
F^{GW}_\g &=& -\g +\frac{1}{2} \log \g + \frac{3}{4}\,,\\
\label{xiGW}
\xi^{GW}_\g &=& \frac{1}{2}(\log \g -\g +1)\,.
\end{eqnarray}
The result (\ref{cstGW}) is shown in \cite{GrossW}. Moreover, (\ref{xiGW}) is  formula (4.14) in \cite{DJ}.
When  $\g < - 1$,  the
 equilibrium measure is 
\begin{equation}\label{GWMeq}
\GW_{\g} (dz) = \frac{|\g|}{\pi}\sin(\tfrac{\theta}{2})\!\  \sqrt{\sin^2(\tfrac{\theta}{2})-\cos^2(\tfrac{\theta_\g}{2})}\!\ 1_{[\pi-\theta_\g, \pi+\theta_\g]}\!\ d\theta\;, (z = e^{\ii \theta} ,\theta\in [0, 2\pi)),\end{equation}
where $\theta_\g$ has the same value as before. It is the same to say that the support of $\GW_\g$ is $\widehat{[\pi-\theta_\g, \pi + \theta_\g]}$.

Let $\widetilde{\mathbb G\mathbb W}_{-\g}\sn$ be the probability measure on $\mathcal H_n$ obtained by pushing forward 
   ${\mathbb G\mathbb W}_{-\g}\sn$  by $\tau$. 
We get 
\[\frac{d\widetilde{\mathbb G\mathbb W}_{-\g}\sn}{d\mathbb Q\sn}(H) := \frac{1}{\widetilde{\mathcal Z}_n(\g)} \exp \left[n\g\,\tr \frac{\mathrm{I}_n-H^2}{\mathrm{I}_n+H^2}\right] [\det (\mathrm{I}_n+H^2)]^{-n} \,.\]
The potential on $\mathbb R$ is
\begin{equation}
\label{potHGW}V_{-\g} (x) = \g \frac{x^2-1}{x^2 + 1} + \log (1+x^2)\,. \end{equation}
For $0 \leq \g \leq 1$ the equilibrium measure (supported by $(-\infty, \infty)$) is
\begin{equation}
\label{tildemuGWs}
\widetilde{\GW}_{-\g} (dx) = \frac{(1 - \g) x^2 + 1+\g}{\pi (x^2+1)^2}\!\ dx
\end{equation}
(for $\g =0$ it is the Cauchy distribution). For $\g > 1$, the equilibrium measure has a compact support:
\begin{equation}
\label{HGWeq1}
\widetilde{\GW}_{-\g} (dx)   
 = \frac{2\sqrt{1+\m^2}}{\pi \m^2}\frac{\sqrt{\m^2-x^2}}{(1+x^2)^2} \ 1_{[-\m , \m]}(x) \!\ dx
\,,\end{equation}
where $\m^2 =(\g -1)^{-1}$.
\subsection{Hua-Pickrell ensemble}
The following distribution has been introduced in \cite{hua1963harmonic} and appears later in \cite{pickrell1987measures}. We also refer to \cite{Ner1}, where the case of a complex parameter is studied.  Further references are \cite{BO} and \cite{BNR}. The Hua-Pickrell ensemble has the following density with respect to the Haar measure on $\mathbb U(n)$: 
\begin{equation}
\label{HP0}\frac{d\mathbb H\mathbb P_{\delta}\sn}{d\mathbb P\sn}(U) := \frac{1}{{\mathcal Z}_n(\delta)}
\left[\det (\mathrm{I}_n-U)\right]^{\bar\delta}
\left[\det (\mathrm{I}_n-\bar U)\right]^{\delta}\,.
\end{equation} 
Here, $\delta$ is a complex parameter such that $\Re \delta > -1/2$. 
Let  $\widetilde{\mathbb H\mathbb P}_{\delta}\sn$ denote the probability measure on $\mathcal H_n$ obtained by pushing forward 
  ${\mathbb H\mathbb P}_{\delta}\sn$ by $\tau$. 
We get 
\[\frac{d\widetilde{\mathbb H\mathbb P}_{\delta}\sn}{d\mathbb Q\sn}(H) := \frac{1}{\widetilde{\mathcal Z}_n(\delta)} \left[\det(\mathrm{I}_n+ H^2)\right]^{-n}\left[\det ( \mathrm{I}_n +iH)\right]^{-\bar \delta} \left[\det ( \mathrm{I}_n -iH)\right]^{- \delta}
\,.\]
A particularly interesting case is the regime $\delta = \d n$, which requires $\Re \d \geq 0$ for integrability. The case $\d= 0$ is of course the same as $\g =0$ in the 
Gross-Witten and corresponds to the Cauchy ensemble.  For simplicity of the computations we will consider here the 
 case $\d > 0$, although it is possible to treat the general case.
 In the framework laid out in Section \ref{sec:RMmodels}, this corresponds to the potential
\begin{eqnarray}
\label{potHP2}
\V_\d(z) = -2\d \log |1-z|\,,
\end{eqnarray}
which is invariant by $z \mapsto \bar z$ and satisfies assumptions (T1) and (T2) and by the remark just after (T3) also this assumption. 
The equilibrium measure is
\begin{equation}
\label{limmeas}
{\HP}_\d (dz) = (1+\d) \frac{\sqrt{\sin^2(\tfrac{\theta}{2}) -
\sin^2(\tfrac{\theta_{\d}}{2})}}{2\pi \!\ \sin(\tfrac{\theta}{2})} \mathbbm{1}_{(\theta_\d, 2\pi-\theta_\d)}(\theta) d\theta, \; (z= e^{\ii \theta}, \theta \in [0, 2\pi])\,,
\end{equation} 
where $
\theta_\d \in (0, \pi)$ is such that 
\begin{equation}\sin \frac{\theta_\d}{2} 
 = \frac{\d}
{ 1 + \d}\,. 
\end{equation}
The support of the equilibrium measure is thus the (symmetric) arc $\widehat{[\theta_\d, 2\pi - \theta_\d]}$.
We have
\begin{eqnarray}
\label{FHP}
F_\d^{HP} &=&(1+\d)^2 \log (1+\d) + \d^2 \log \d \\
&&-\frac{1}{2}(1+2\d)^2 \log (1+ 2\d) + 2\d^2 \log 2\,,\notag \\
\label{xiHP}
\xi_\d^{HP} &=& (1+\d) \log(1+\d) - \frac{1+2\d}{2}\log (1+ 2\d)\,.
\end{eqnarray}
The orthogonal polynomials are the Geronimus polynomials with constant Verblunsky coefficients
\begin{equation}
\alpha_k \equiv \gamma_\d, \qquad k \geq 0\,, 
\end{equation}
where
\begin{equation}
\label{defgammad}
\gamma_\d := - \frac{\d}{1+\d}, \qquad k \geq 0\,.
\end{equation}
Pushing forward this measure on the set $\mathcal H_n$ of $n\times n$ Hermitian matrices, we get
 the potential
\begin{equation}
\label{potHP3}
 V_\d(x) = (1+ \d) \log (1+x^2)\,.
\end{equation}
This model is sometimes called the modified Cauchy ensemble, see \cite{Forrester}, \cite{Jo}, \cite{Mizoguchi}, \cite{pastur_shcherbina} (Problem 11.4.15), or the  Lorentzian ensemble \cite{brouwer1995generalized}. The equilibrium measure
on the real line is 
\begin{equation}
\label{param21}
\widetilde{\HP}_\d(dx) 
= \frac{1}{\pi(\sqrt{1+ \p^2} -1)}\frac{\sqrt{\p^2-x^2}}{1+x^2} \ 1_{[-\p , \p]}(x) \!\ dx\,,
\end{equation}
where  $\p^2  = (1 + 2\d) \d^{-2}$ (see \cite{Blower} Prop. 11.2.2, p. 359).
 Moreover
\begin{eqnarray}
\label{Fcauchy1}
\widetilde{F}_\d^{HP} &=& (1+\d)^2 \log (1+ \d) + \d^2 \log \d \\
&&-\frac{1}{2}(1+2\d)^2 \log (1+ 2\d) + (2\d^2-1) \log 2\,,\notag \\
\label{xicauchyp1}
\widetilde{\xi}_\d^{HP} &=& \left(\d +\frac{1}{2}\right)\log (1+ 2\d) -\d\log \d -(1+ 2\d)\log 2\,.
\end{eqnarray}
\begin{rem}
The corresponding Jacobi coefficients of the tridiagonal representation are
\begin{eqnarray}
a_1 = \sqrt{\frac{2(1+ 2\d)}{(1+ \d)^3}} &,& \ a_k = \frac{1+ 2\d}{(1+ \d)^2}\ \ (k > 1)\,,\\
b_1 = -\frac{2\d}{1+ \d} &,& \ b_k = -2 \frac{\d^2}{(1+ \d)^2}\ \ (k > 1)\,.
\end{eqnarray}
 We did not find the corresponding values in the literature.
\end{rem}
As an application of Corollary \ref{LDPmuuU} and Proposition \ref{LDPmuuR}, we have the following result, collecting all the LDPs for the empirical spectral measure as in \eqref{empiricallaw} in our basic models. 
\begin{cor} 
\label{ESDgene}
\begin{enumerate}
\item For any $\g \in \mathbb R$, 
 the sequence of distributions of $(\mu_\u\sn)$ under $\mathbb G\mathbb W\sn_\g$ satisfies the LDP in ${\mathcal M}_1(\mathbb T)$, with speed $n^2$ and good rate funtion $I_{\V}$ with $\V = \V_\g$ given by (\ref{potGW}).

\item For any $\d >0$,  the sequence of distributions of $(\mu_\u\sn)$ under $\mathbb H\mathbb P_{\d n}\sn$ satisfies the LDP in ${\mathcal M}_1(\mathbb T)$, with speed $n^2$ and good rate funtion $I_{\V}$ with $\V = \V_\d$ given by (\ref{potHP2}).

\item For any $\g \in \mathbb R$, 
the sequence of distributions of $(\mu_\u\sn)$ under $\widetilde{\mathbb G\mathbb W}\sn_\g$ satisfies the LDP in ${\mathcal M}_1(\mathbb R)$, with speed $n^2$ and good rate funtion $I_{V}$ with $V = V_{-\g}$ given by (\ref{potHGW}).
\item For any $\d >0$,  the sequence of distributions of $(\mu_\u\sn)$ under $\widetilde{\mathbb H\mathbb P}_{\d n}\sn$ satisfies the LDP in ${\mathcal M}_1(\mathbb R)$, with speed $n^2$ and good rate funtion $I_{\V}$ with $V = V_\d$ given by (\ref{potHP3}).
\end{enumerate}
\end{cor}

Point 1. is in \cite{HiaiP} p. 225 and point 2. is in \cite{BNR} Theorem 5.5. The points 3. and 4. are obtained carrying the results to the real line by the Cayley transform.
\section{LDP for spectral measures}
\label{sLDDPP}
\subsection{Measure encoding approach}
In this subsection, we state LDPs for the weighted measures given in \eqref{spectralmeasure}. 
They are elements of $\mathcal M_1 (\mathbb T)$.
We first recall the main theorem of \cite{GaNaRo} on  $\mathbb R$, then we state the LDP on $\mathbb T$ improving the result on $\mathbb R$ with weaker assumptions. To begin with,  we recall the definition of the Kullback-Leibler divergence, with a slight generalization for sub-probabilities.

Let $\mu$ be a probability measure and $\nu$ be a non-zero sub-probability measures on some measurable space. The Kullback-Leibler divergence between $\mu$ and $\nu$ is given by
\begin{equation}
\label{kukul}
\mathcal{K}(\mu|\nu) = \int \log \left( \frac{d \mu}{d\nu} \right) d\mu 
\end{equation}
if $\mu$ is absolutely continuous with respect to $\nu$ and $\log\frac{d \mu}{d\nu}\in L^1(\mu)$. Further, set $\mathcal{K}(\mu|\nu)=\infty$ otherwise. In our LDP, the rate function will involve the reversed Kullback-Leibler distance, where $\mu$ will be the reference measure and $\nu$ is the argument.
Recall the definition of the set $\Sr_1^\mathbb{R}(\am,\ap)$ given in the introduction. It consists in  probability measures 
\begin{align}\label{muinS0R}
\mu = \mu_{|I} +  \sum_{i=1}^{N^+} \gamma_i^+ \delta_{\lambda_i^+} + \sum_{i=1}^{N^-} \gamma_i^- \delta_{\lambda_i^-}\,.
\end{align}
In our extension of the Killip-Simon sum rule we will also consider reference measures supported by the whole real line. To keep coherent notations, we write $\Sr_1^\mathbb{R}(-\infty,\infty)$ for the set of probability measures with support $\mathbb{R}$. In this case, $N^+=N^-=0$. In the same vein, we define $S_{\leq 1}^{\mathbb R}(\alpha^-, \alpha^+)$ in the case of subprobabilities. Notice that this last set may be seen as $\mathcal S_1^{\bar{\mathbb R}} (\alpha^- , \alpha^+)$.

We now introduce the analogous framework on the circle. 
If $[\alpha^- , \alpha^+]$ is an interval in $(0, 2\pi)$, let 
$I=\widehat{[\alpha^- , \alpha^+]}$
and let $\Sr_1^{\mathbb T} = \Sr_1^{\mathbb T}(\am,\ap)$ be the set of all probability measures $\mu$ on $\mathbb T$ with 
\begin{itemize}
\item[(i)] $\operatorname{supp}(\mu) = J \cup \{e^{\ii \theta_i^-}\}_{i=1}^{N^-} \cup \{e^{\ii \theta_i^+}\}_{i=1}^{N^+}$, where $J\subset I$, $N^-,N^+\in\N_0\cup\{\infty\}$ and $\theta_i^\pm \in [0,2\pi)$. Furthermore,
\begin{align*}
0\leq\theta_1^-<\theta_2^-<\dots <\am \quad \text{and} \quad \theta_1^+>\theta_2^+>\dots >\ap .
\end{align*}
\item[(ii)] If $N^-$ (resp. $N^+$) is infinite, then $\theta_j^-$ converges towards $\am$ (resp. $\theta_j^+$ converges to $\ap$).
\end{itemize}
We will also write $\lambda_i^\pm =e^{\ii \theta^\pm_i}$ as in the real case for the outlying support points. For a measure $\mu\in \Sr_1^{\mathbb T}(\am,\ap)$ we may write it as in \eqref{muinS0R}. 
%
Like in the real case, we write $\Sr_1^{\mathbb T}(0,2\pi)$ for the probability measures supported by $\mathbb T$. 
It should be clear that the Cayley transform carries  $\{\mu \in \Sr_1^{\mathbb T}| \mu(1)=0\}$ onto $\Sr_1^{\mathbb R}$ and $\mathcal S_1^{\mathbb T}$ onto $\mathcal S_{\leq 1}^{\mathbb R}$. 
Furthermore, as the circle is rotationally invariant, classifying an outlier in $(\theta_i^+)$ or $(\theta_i^-)$ is essentially arbitrary. Nevertheless,  it is consistent with our measure mapping. 
We endow the sets $\Sr_1^{\mathbb T}$ and $\Sr_1^{\mathbb R}$ with the weak topology and $\Sr_{\leq1}^{\mathbb T}$ and $\Sr_{\leq1}^{\mathbb R}$ with the vague topology and the corresponding Borel $\sigma$-algebra.

We need 
one more definition 
 in order to formulate the general result. 
Recall that  $\mathcal{J}_V$ has been defined in  assumption (A3). We define, in the general case, the rate function for the extreme eigenvalues,
\begin{align}
\label{rate0}
\mathcal{F}_V^+(x) & = \begin{cases}
\mathcal{J}_V(x) - \inf_{\xi \in \R} \mathcal{J}_V(\xi) & \text{ if } x\geq \alpha^+, \\
\infty & \text{ otherwise, } 
\end{cases} \\
\mathcal{F}_V^-(x) & = \begin{cases}
\mathcal{J}_V(x) - \inf_{\xi \in \R} \mathcal{J}_V(\xi) & \text{ if }  x \leq \alpha^-, \\
\infty & \text{ otherwise. } 
\end{cases}
\end{align}
On the unit circle, we have similar notations, with $V$ replaced by $\V$.
Notice that if $\mathcal V(1) < \infty$, then $\mathcal F_{\mathcal V}(1) < \infty$. 
In this case let us denote
\begin{equation}
\label{kappa}
\kappa_{\mathcal V} = \mathcal F_{\mathcal V}(1) .
\end{equation}
\begin{prop}[\cite{GaNaRo} Thm. 3.1] \label{MAINR}
Assume that the potential $V$ satisfies assumptions (R1s), (R2) and (R3). Then the sequence of spectral measures $\mun$ under $\Pnv$ satisfies the LDP with speed $n$ and good rate function
\begin{align*}
\mathcal{I}_V(\mu) = \mathcal{K}(\mu_V\!\ |\!\ \mu) + 
\sum_{i=1}^{N^+} {\mathcal F}_{V}^+ (\lambda_i^+)  +  \sum_{i=1}^{N^-} {\mathcal F}_{V}^- (\lambda_i^-)
\end{align*}
if $\mu \in \mathcal{S}^{\mathbb R}_1(\alpha^-,\alpha^+)$ and $\mathcal{I}_V(\mu) = \infty$ otherwise. 
\end{prop}
On the unit circle, we claim:
\newpage
%
\begin{thm} 
\label{MAINT} \ 
\begin{enumerate}
\item
Assume that the potential $\mathcal V$ satisfies (T1) and that the support of $\mu_{\mathcal V}$ is $\mathbb T$. Then the sequence of spectral measures $\mun$ under $\mathbb P\sn_{\V}$ satisfies the LDP in $\mathcal M_1(\mathbb T)$ with speed $n$ and good rate function
\begin{equation}
\mathcal I_{\mathcal V}( \mu) = \mathcal{K}(\mu_\V\!\ |\!\ \mu)\,.
\end{equation}
\item
Assume that the potential $\V$ satisfies the assumptions (T1), (T2) and (T3). Then, the sequence of spectral measures $\mun$ under $\mathbb P\sn_{\V}$ satisfies the LDP in $\mathcal M_1(\mathbb T)$ with speed $n$ and good rate function
\begin{align}
\label{MAINTf}\mathcal{I}_\V(\mu) = \mathcal{K}(\mu_\V\!\ |\!\ \mu) + 
\sum_{i=1}^{N^+} {\mathcal F}_{\V}^+ (\lambda_i^+)  +  \sum_{i=1}^{N^-} {\mathcal F}_{\V}^- (\lambda_i^-)
\end{align}
if $\mu \in \mathcal{S}^{\mathbb T}_1(\alpha^-,\alpha^+)$ and $\mathcal{I}_\V(\mu) = \infty$ otherwise. 
\end{enumerate}
\end{thm}
To transfer the LDP in Theorem \ref{MAINT} to the real line we use  the 
mapping $\hat\tau$ given in (\ref{hattau}). 
We get the following corollary.
\begin{cor}\label{cor:weakassump} \ 
\begin{enumerate}
\item Assume that the potential $V$ satisfies the assumption (R1w) and that the support of $\mu_V$ is $\mathbb R$. Then, the sequence of spectral measures $\mun$ under $\Pnv$ satisfies the LDP in $\mathcal M_{\leq 1} (\mathbb R)$ with speed $n$ and good rate function
\begin{equation}
\mathcal{I}_V(\mu) = \mathcal{K}(\mu_V\!\ |\!\ \mu)\,. 
\end{equation}
\item Assume that the potential $V$ satisfies the assumptions (R1w), (R2) and (R3). 
Then the sequence of spectral measures $\mun$ under $\Pnv$ satisfies the LDP  in $\mathcal M_{\leq 1}(\mathbb R)$ with speed $n$ and good rate function
\begin{align*}
\mathcal{I}_V(\mu) = \mathcal{K}(\mu_V\!\ |\!\ \mu) + 
\sum_{i=1}^{N^+} {\mathcal F}_{V}^+ (\lambda_i^+)  +  \sum_{i=1}^{N^-} {\mathcal F}_{V}^- (\lambda_i^-) + \kappa_{\mathcal V}\mathbbm{1}_{  \mu(\mathbb R) < 1}
\end{align*}
if $\mu \in \mathcal{S}^{\mathbb R}_{\leq 1}(\alpha^-,\alpha^+)$ and $\mathcal{I}_V(\mu) = \infty$ otherwise. 
\end{enumerate}
\end{cor}

\proof We only prove the second point, since the other one is more straighforward. 
Under $\mathbb P_{\mathcal V}^{(n)}$, we consider the two random measures
\begin{eqnarray*}
\nu^{(n)} = \sum_{k=1}^n  \w_k \delta_{\zeta_k} \in \mathcal{M}_1(\mathbb T)\ \ \hbox{and} \ \  
\mu^{(n)} = \hat \tau( \nu^{(n)}) = \sum_{k=1}^n  \w_k \delta_{\tau(\zeta_k)}\in \mathcal{M}_{\leq 1}(\mathbb R)\,.
\end{eqnarray*}The mapping $\hat\tau$ is continuous, and $\mathcal I_{\mathcal V}$ is good. We may apply the contraction principle (Theorem 4.2.1 in \cite{demboz98}). We obtain the LDP in $\mathcal M_{\leq 1}(\mathbb R)$ with good rate function
\[\widehat{\mathcal I}(\mu) = \inf_{\nu : \hat\tau(\nu) = \mu} \mathcal I_{\mathcal V} (\nu) \,. \] 
Actually only those $\nu$ such that $\mathcal I_{\mathcal V}(\nu)$ is finite contribute to the infimum. Therefore,
$\nu \in \mathcal S_1^{\mathbb T}(\alpha_ {\mathcal V}^-, \alpha_{\mathcal V}^+)$ implies that $\hat\tau (\nu) \in \mathcal S_{\leq 1}^{\mathbb R}( \alpha_V^-, \alpha_V^+)$ with $\alpha_V^\pm = \tau(\alpha_{\mathcal V}^\pm)$. Under our assumptions, $\mu_{\mathcal V}$ has no atom at $1$ and $\mu_V = \hat\tau (\mu_{\mathcal V}).$ For a $\nu$ such as above, we have by pushforward by $\hat\tau$ 
\[\mathcal K( \mu_{\mathcal V} | \nu) = \mathcal K(\hat\tau (\mu_{\mathcal V}) | \hat\tau(\nu))  =  \mathcal K( \mu_{V} | \mu)\,.\]
Moreover, the outliers of $\nu$ different from $1$  are carried to outliers of $\mu$, and $\mathcal F_{\mathcal V}(\zeta) = \mathcal F_V(\tau(\zeta))$. 
Besides, when $\nu$ has an outlier at $1$, say $\nu = \nu_0 + a \delta_{1}$ then $\mu(\mathbb R) = \nu_0(\mathbb T) = 1-a$, and the contribution of $1$ in $\mathcal I_{\mathcal V}(\nu)$ hence in $\widehat{\mathcal I}(\mu)$ is  $\kappa_{\mathcal V}$.  This proves that $\widehat{\mathcal I} = \mathcal I_V$ and ends the proof of the corollary. $\Box$

As a consequence, we have for our models the following results.

\begin{cor}   \ 
\label{appliHP} 
\begin{enumerate} 
\item
Under $\mathbb H\mathbb P_{n\d}\sn$, the sequence of spectral measures $\mu\sn$ satisfies the LDP 
in $\mathcal M_1(\mathbb T)$
with speed $n$ and good rate function $\mathcal I_\V$
where $\mu_\V = \HP_\d$ is given in (\ref{limmeas}) and $\mathcal{F}^\pm_\V = \mathcal{F}^\pm_{HP}$, where for $0 < \theta \leq \theta_\d$
\begin{eqnarray}
\mathcal F_{HP}^- (e^{\ii \theta}) := \int _\theta^{\theta_\d} (1+\d) \frac{\sqrt{\sin^2\big(\theta_\d/2\big) -
\sin^2(\varphi/2)}}{2\sin(\varphi/2)}\!\   d\varphi
\end{eqnarray}
and for $\theta \in (2\pi - \theta_\d, 2\pi)$, $\mathcal F_{HP}^+ (e^{\ii \theta})  := \mathcal F_{HP}^- (e^{-\ii \theta})$. 
\item Under $\widetilde{\mathbb H\mathbb P}_{n\d}\sn$, the sequence of spectral measures $\mu\sn$ satisfies the LDP with speed $n$ and good rate function $\mathcal I_V$, where $\mu_V = \widetilde{\HP}_\d$ is given by (\ref{param21}) and $\mathcal{F}^\pm_V = \mathcal{F}^\pm_{HP}$, where for $x \geq \p$
\begin{eqnarray}
\mathcal F_{\widetilde{HP}}^+ (x) = \int^x_\p  \frac{2}{\sqrt{1+\p^2} -1}\frac{\sqrt{\xi^2 -\p^2}}{1+\xi^2} \  \!\ d\xi\,,
\end{eqnarray}
and $\mathcal F^-_{\widetilde{HP}} (x) = \mathcal F^+_{\widetilde{HP}} (-x) $ for $x \leq -\p$.
\end{enumerate}
\end{cor}
\begin{cor} \ 
\label{appliGW} 
\begin{enumerate} 
\item
Under $\mathbb G\mathbb W\sn_\g$, $\g \leq 1$, the sequence of spectral measures $\mu\sn$ satisfies the LDP  in $\mathcal M_1(\mathbb T)$ with speed $n$ and good rate function $\mathcal{I}_\V$ where $\mu_\V = \GW_\g$ is given in \eqref{GW-} and \eqref{GWMeq}.  
\\ If $|\g| \leq 1$ there is no outlier and the rate function reduces to
\[\mathcal{I}_\V(\mu) = \mathcal K(\GW_\g | \mu)\,.\]
If $\g < -1$, we have $\mathcal F^\pm_\V= \mathcal F^\pm_{GW}$, where for $0 < \theta < \pi- \theta_\g$
\[\mathcal F_{GW}^- (e^{\ii \theta}) = \int_\theta^{\pi-\theta_\g} 2|\g| \sin\frac{\varphi}{2}\!\  \sqrt{\cos^2\frac{\theta_\g}{2} - \sin^2\frac{\varphi}{2}}\!\  d\varphi   = 4\int^{\sqrt{  |\g|} \cos \frac{\theta}{2}}_1 \sqrt{u^2 - 1} \!\ du  \,,\]
and 
$\mathcal F_{GW}^+ (e^{\ii \theta}) = \mathcal F_{GW}^- (e^{-\ii \theta}) $ if $\pi + \theta_\g < \theta < 2\pi$. 
\item 
Under $\widetilde{\mathbb G\mathbb W}\sn_{-\g}$, $\g\geq 0$, the sequence of spectral measures $\mu\sn$ satisfies the LDP  with speed $n$ and good rate function $\mathcal{I}_V$ with $\mu_V = \widetilde\GW_{-\g}$ as in \eqref{tildemuGWs} and \eqref{HGWeq1}.
\\  If $0\leq \g \leq 1$, the support of $\mu_V$ is $\mathbb{R}$, the LDP is in $\mathcal M_{\leq 1}(\mathbb R)$ and the rate function is
\[\mathcal{I}_V(\mu)=  \mathcal K(\widetilde{\GW}_{-\g} | \mu)\,.\] 
If $\g > 1$, the LDP is in $\mathcal M_{\leq 1}(\mathbb R)$. We have $\mathcal F^\pm_V= \mathcal F^\pm_{\widetilde{GW}}$ where for $ x > \m$
\[\mathcal F_{\widetilde{GW}}^+ (x) = \int_\m^x 
\frac{4\sqrt{1+\m^2}}{\m^2}\frac{\sqrt{\xi^2 -\m^2}}{(1+\xi^2)^2}  \!\ d\xi = 4 \int_\m^{\frac {x|\g|}{\sqrt{1+x^2}}}\sqrt{u^2 - 1} \!\ du  \,,\]
and for $x < -\m$, $\mathcal F_{\widetilde{GW}}^- = \mathcal F_{\widetilde{GW}}^+ (-x) $.
\end{enumerate}
\end{cor}

\subsection{Verblunsky coefficient encoding approach}
\label{suveve}
To begin with, let us recall the simplest example.
It is the Circular Unitary Ensemble where $\mathbb U(n)$ is equipped with the Haar measure. Then 
the Verblunsky coefficients 
are independent. 
More precisely, Killip and Nenciu proved in \cite{Killip1} that 
the $n$-tuple $\alpha\sn :=\left(\alpha_0, \dots, \alpha_{n-1}=e^{\ii\phi}\right)$ has the distribution
\begin{equation}
\label{ref}dP_0\sn (\alpha_0, \dots, \alpha_{n-1}) = \left(\otimes_{k=0}^{n-2} \eta_{n-k-2}(d\alpha_k)\right)\otimes\frac{d\phi}{2\pi}\end{equation}
where, for $r > -1$
\begin{equation}
\label{defeta}
\eta_r(d\alpha) := \frac{r+1}{\pi}\left(1 - |\alpha|^2\right)^r \, d\alpha\,,
\end{equation}
and $d\alpha$ is the Lebesgue measure on the unit disk. 
From that, it is deduced in \cite{gamboacanonical}, Section 5.2 that the family of distributions of   $\mu\sn$ under CUE$(n)$ satisfies the LDP (in $\mathcal M_1(\mathbb T)$ equipped with the weak topology) with speed $n$ and good rate function
\[I^0 (\mu) = \sum_{k=0}^\infty -\log (1-|\alpha_k|^2)\,,\]
when $\alpha_k, k\geq 0$ are the Verblunsky coefficients of $\mu$. In the Hua-Pickrell case, the Verblunsky coefficients are no more independent (except when $\d = 0$). To recover a structure of independence, it is necessary to introduce the so-called {\it deformed Verblunsky coeffficients}. Given a measure $\mu \in \mathcal M_1(\mathbb T)$ with at least $n$ distinct support points and monic orthogonal polynomials $\phi_0,\dots ,\phi_{n-1}$, define 
\begin{equation}
\label{phioverphi}
b_k = \frac{\phi_k(1)}{\phi_k^*(1)}\quad \hbox{and} \quad  \gamma_k = \bar\alpha_k (b_k)^{-1}\, , \quad k =0, \dots, n-1\,.\end{equation}
This is equivalent to the recursive definition
\begin{equation}
\label{mVc}
\gamma_0 = \bar\alpha_0 , \quad  \gamma_k = \bar\alpha_k \prod_{j=0}^{k-1} \frac{1 - \bar\gamma_j}{1 - \gamma_j}\, , \quad k = 1, \dots, n-1\,.
\end{equation}
For a more detailed description and meaning of these quantities we refer to 
\cite{BNR}, Section 2.2. In Theorem 3.2 therein, it is proved that under $\mathbb H \mathbb P_\delta^{(n)}$, the random variables $\gamma_0\sn, \dots, \gamma_{n-1}\sn$ are independent and 
for $k = 0, \dots, n-2$, the density of $\gamma_k\sn$ on $\mathbb D$ is
\begin{equation}
\label{BNR1}
\frac{\Gamma(n-k+ \delta)\Gamma(n-k + \bar\delta)}{\pi \Gamma(n-k-1)\Gamma(n-k + \delta + \bar\delta)} (1 - |z|^2)^{n-k-2} (1-z)^{\bar\delta}(1-\bar z)^{\delta}\,,
\end{equation}
and $\gamma_{n-1}\sn \in \mathbb{T}$ has the density 
\label{BNR2}\begin{equation}
\frac{\Gamma(1 + \delta)\Gamma(1 + \bar\delta)}{\Gamma(1 + \delta + \bar\delta)}( 1 - \zeta)^{\bar\delta}(1- \bar\zeta)^{\delta}\end{equation}
with respect to the Haar measure on $\mathbb T$.

When $\delta = n \d$, $\d \geq 0$, a straightforward study of the density  (\ref{BNR1}) leads to a LDP for $\gamma_j^{(n)}$ for $j$ fixed. It is a particular case of the matricial result (Proposition \ref{lemmat}) proved in Section \ref{sec:proofoflemmat}. 

\begin{lem}
\label{LDPalphaHP}
For fixed $k$, $(\gamma_0\sn, \gamma_1\sn, \dots, \gamma_k\sn)_{n \geq k}$ satisfies under $\mathbb H\mathbb P_{n\d}\sn$ the LDP in  $\bar{\mathbb D}^k$ with speed $n$ and good rate function
\[I_k(\gamma_0, \dots, \gamma_k) = \sum_{j=0}^k H_\d(\gamma_j)\,,\]
where
\begin{align}
H_\d(\gamma) &= -\log (1 - |\gamma|^2)  - 2 \d \log |1- \gamma|  + H_\d(0), \\
\label{H1}
H_\d(0) &= \log (1 - \gamma_\d^2)  - 2 \d \log (1- \gamma_\d)
= (1 + 2\d) \log (1+2\d) - 2(1+\d) \log(1+\d)\,.
\end{align}
\end{lem}

Note that $H_\d$ has its unique minimum at $\gamma_\d$ 
, the (deformed) Verblunsky coefficient of the Hua-Pickrell distribution (\ref{defgammad}). 
Using the classical method of projective limits (see the proof in Section \ref{sec:proofs}), this allows to claim:

\begin{thm}
\label{JdHP}
Under  $\mathbb H\mathbb P_{n\d}\sn$, the sequence of measures $\mu\sn$ satisfies the LDP in $\mathcal M_1(\mathbb T)$ with speed $n$ and good rate function
\[J_\d^{HP}(\mu) = \sum_{k=0}^\infty  H_\d(\gamma_k)\,.\]
if $\mu$ is non-trivial and infinite otherwise. 
\end{thm}

Of course, if $\mu$ has only $n$ support points, only the first $n$ (deformed) Verblunsky coefficients can be defined. Then $\alpha_{n-1} \in \mathbb{T}$ and also $\gamma_{n-1} \in \mathbb{T}$, which implies $H_\d(\gamma_{n-1})=\infty$ and the rate function is infinite. 

In the Gross-Witten case, the Verblunsky coefficients are not independent (except when $\g = 0$).  More precisely, the joint distribution is given by the following lemma.

\begin{lem}
\label{lemdGW}
The law of $\left(\alpha_0\sn, \dots, \alpha_{n-1}\sn\right)$ under $\mathbb G\mathbb W_\g\sn$ 
 is given by 
\begin{equation}
\label{densalpha}dP_\g\sn\left(\alpha_0, \dots, \alpha_{n-1}=e^{i\phi}\right) = {{\mathcal Z}_n(\g)}^{-1} \exp \left[n\g \Re \left(\alpha_0 - \sum_{k=0}^{n-1} \alpha_k \bar \alpha_{k-1}\right) \right] \left(\otimes_{k=0}^{n-2} \eta_{n-k-2}(d\alpha_k)\right)\otimes\frac{d\phi}{2\pi}\,.\end{equation}
\end{lem}

\proof By definition, we have
\[\frac{d\mathbb G\mathbb W_{\g}\sn}{d\mathbb P\sn}(U)= {{\mathcal Z}_n(\g)}^{-1} \exp \left(\frac{n\g}{2} \tr\!\ (U + U^\dagger)\right) = {{\mathcal Z}_n(\g)}^{-1} \exp \left(n\g \Re(\tr\!\ U)\right)\,.\]
Now, from  the CMV representation (\ref{favardfini}),  we get
\begin{equation}
\label{trace}\tr\!\ U = -\alpha_0 + \sum_{k=0}^{n-1} \alpha_k \bar\alpha_{k-1}\,,\end{equation}
(see also  Simon \cite{simon05} p. 273). It remains to use (\ref{ref}). \hfill $\Box$

Given the explicit density in Lemma \ref{lemdGW}, we may conjecture a LDP for the spectral measure in terms of its Verblunsky coefficients. B. Simon (personal communication) has notified us of a forthcoming paper (\cite{BSZ1}) with J. Breuer and O. Zeitouni in which LDPs for certain  ensembles on the unit circle, and in particular for  the Gross-Witten example, are obtained (see also Section 6 of \cite{BSZ}). 

\begin{conj}\label{GWconjLDP}
Under  $\mathbb G\mathbb W_{\g}\sn$, the sequence of measures $\mu\sn$ satisfies the LDP in $\mathcal M_1(\mathbb T)$ with speed $n$ and rate function
\begin{align*}
J_\g^{GW}(\mu) = H(\g) - \g \Re \left(\alpha_0 - \sum_{k=1}^\infty \alpha_k\bar\alpha_{k-1}\right) - \sum_{k=0}^\infty \log (1 - |\alpha_k|^2).
\end{align*}

\end{conj}

\section{Sum rules from large deviations}
\label{sumrule}

\subsection{Hua-Pickrell case}

Our new sum rule is a straightforward consequence of Theorem \ref{JdHP} and Corollary \ref{appliHP}.

\begin{thm}
\label{sumruleHP}
Let $\mu \in \mathcal M_1(\mathbb T)$ with infinite support and let  
$(\gamma_k)_{k \geq 0}\in {\mathbb{D}}^{\mathbb N}$ be the sequence of its deformed Verblunsky coefficients. Then, for any $\d \geq 0$, we have $\sum_{k=0}^\infty H_\d(\gamma_k)=\infty$ if $\mu \notin \Sr_1^{\mathbb T}(\theta_\d,2\pi-\theta_\d)$. If $\mu \in \Sr_1^{\mathbb T}(\theta_\d,2\pi-\theta_\d)$, we have
\ben
 \mathcal K(\HP_\d | \mu) + \sum_{i=1}^{N^+} \mathcal F_{HP}^+(\lambda_i^+) + \sum_{i=1}^{N^-} \mathcal F_{HP}^-(\lambda_i^-) = \sum_{k=0}^\infty H_\d(\gamma_k)\,,
\een
where both sides may be infinite simultaneously.
\end{thm}

As we wrote in the introduction, an essential consequence of a sum rule are \emph{gems}, equivalent conditions for finiteness of the rate function. For the Hua-Pickrell case, Theorem \ref{sumruleHP} and an expansion of $H_\d$ in the neighbourhood of $\gamma_d$ gives the following corollary.

\begin{cor}
\label{gemHP}
Let $\mu$ be a probability measure on $\mathbb{T}$ with infinite support and deformed Verblunsky coefficients $(\gamma_k)_{k \geq 0}\in {\mathbb{D}}^{\mathbb N}$. Then
\begin{align*}
\sum_{k=1}^\infty |\gamma_k - \gamma_\d 
|^2 < \infty
\end{align*}
(that is, $J_{\d}^{HP}(\mu)< \infty$) if and only if
\begin{enumerate}
\item 
  $\mu \in \Sr_1^{\mathbb T}(\theta_\d,2\pi-\theta_\d)$
\item $\sum_{i=1}^{N^+} (\theta_i^+ - 2\pi+\theta_\d)^{3/2} + \sum_{i=1}^{N^-} (\theta_\d - \theta_i^- )^{3/2}  < \infty$ and if $N^->0$, then $\theta_1^- > 0$. 
\item If $d\mu(\theta) = f(\theta)\frac{d\theta}{2\pi}+d\mu_s(\theta)$ is the decomposition of $\mu$ with respect to the Lebesgue measure, then
\begin{align*}
\int_{\theta_\d}^{2\pi-\theta_\d}  \frac{\sqrt{\sin^2(\tfrac{\theta}{2}) -
\sin^2(\tfrac{\theta_{\d}}{2})}}{2\pi \!\ \sin(\tfrac{\theta}{2})} \log(f(\theta)) d\theta >-\infty .
\end{align*}
\end{enumerate}

\end{cor}

\textbf{Proof of Corollary \ref{gemHP}:}  Point \emph{1} to \emph{3} are equivalent to finiteness of the left side of the equation in Theorem \ref{sumruleHP}. Indeed, we have 
\begin{align*}
\mathcal F_{HP}^- (e^{\ii \theta}) = \int _\theta^{\theta_\d} (1+\d) \frac{\sqrt{\sin^2\big(\theta_\d/2\big) -
\sin^2(\varphi/2)}}{2\sin(\varphi/2)}\!\   d\varphi = c_\d(\theta_\d-\theta)^{3/2} + o((\theta_\d-\theta)^{3/2} )
\end{align*}
as $\theta \to \theta_\d$, so 
the second point is equivalent to $\sum_{i=1}^{N^+} \mathcal F_{HP}^+(\lambda_i^+) + \sum_{i=1}^{N^-} \mathcal F_{HP}^-(\lambda_i^-)$ being finite.
The third point is equivalent to $\mathcal K( \HP_\d | \mu)$ being finite. Corollary \ref{gemHP} follows then from the equality in Theorem \ref{sumruleHP} since an  expansion of $H_\d$ in the neighbourhood of $\gamma_d$ gives 
\[\frac{(1+\d)^3}{(1+2\d)^2} |h|^2 + o(|h|^2) \leq H_\d (\gamma_\d + h) \leq \frac{(1+\d)^3}{1+2\d} |h|^2 + o(|h|^2)\,.\] \hfill $\Box$

\subsection{Gross-Witten case}

As we saw above, we do not have independence of 
 Verblunsky coefficients and could not succeeded in finding a LDP directly  with this encoding. Nevertheless, Simon gave in this frame a sum rule (\cite{simon05}). Here, the reference measure is supported by the full circle $\mathbb{T}$ and there is no contribution of outliers.

\begin{prop}[\cite{simon05} Thm. 2.8.1]
\label{5.3}
Let $\mu\in \mathcal M_1 (\mathbb T)$ 
 with Verblunsky coefficients $(\alpha_k)_{k \geq 0}\in {\mathbb{D}}^{\mathbb N}$. Then
\begin{equation}
\label{sumrulegw1}
\mathcal K(\GW_{-\mathbf 1}| \mu) = 1 - \log 2 +\Re (\alpha_0) + \frac{|\alpha_0|^2}{2} + \frac{1}{2}\sum_{k=1}^\infty |\alpha_k - \alpha_{k-1}|^2+ \sum_{k=0}^\infty h(\alpha_k)\,,
\end{equation}
where 
\[h(\alpha) = - \log (1 - |\alpha|^2) - |\alpha|^2\,.\]
In particular, 
\begin{equation}
\label{gemgw1}
\mathcal K(\GW_{-\mathbf 1}| \mu) <\infty \Longleftrightarrow 
\sum_{k=0}^\infty |\alpha_{k+1} - \alpha_k|^2 +|\alpha_k|^4 < \infty\,.\end{equation}
\end{prop}

As an easy corollary, we have

\begin{cor}
\label{coreasy}
Let $\mu$ be a probability measure on $\mathbb{T}$ with Verblunsky coefficients $(\alpha_k)_{k \geq 0}\in \mathbb D^{\mathbb N}$. Then,
for $0 \leq\g < 1$, 
we have
\begin{align}
\label{sumrulegwg}
\mathcal K(\GW_{-\g}| \mu) &=  H(\g) +\g\left(\Re (\alpha_0) + \frac{|\alpha_0|^2}{2} + \frac{1}{2}\sum_{k=1}^\infty |\alpha_k - \alpha_{k-1}|^2\right)\\
 &\quad \notag + \sum_{k=0}^\infty - \log (1- |\alpha_k|^2) - \g |\alpha_k|^2\,,
\end{align}
where 
\begin{equation}
\label{defH}H(\g) := 
 1 -\sqrt{1-\g^2} +  \log \frac{1+\sqrt{1-\g^2}}{2} 
\,.\end{equation}
In particular, we have
\begin{equation}
\label{gemgwg}
\mathcal K(\GW_{-\g}| \mu) <\infty \Longleftrightarrow 
\sum_{k=0}^\infty |\alpha_k|^2 < \infty\, ,\end{equation}
which also follows from the Szeg\H{o}-Verblunsky sum rule.
\end{cor}

\begin{rem} In the way to prove (\ref{sumrulegw1}), Simon shows the equivalent relation:
\begin{equation}
\label{fortrace1}
\mathcal K(\GW_{-\mathbf 1}| \mu) = 1 - \log 2 +\Re \left( \alpha_0- \sum_{k=1}^\infty \alpha_k \bar\alpha_{k-1}\right)  + \sum_{k=0}^\infty - \log (1 - |\alpha_k|^2)\,.
\end{equation}
In the same vein, (\ref{sumrulegwg}) is equivalent to
\begin{equation}
\label{fortraceg}
\mathcal K(\GW_{-\g} | \mu) = H(\g)+\g \Re \left( \alpha_0- \sum_{k=1}^\infty \alpha_k \bar\alpha_{k-1}\right)+ \sum_{k=0}^\infty - \log (1 - |\alpha_k|^2)\,.
\end{equation}
\end{rem}

For $|\g| > 1$, we may still conjecture a sum rule. The left hand side of such an identity would be given by the rate function of the LDP for the spectral measure encoded by the eigenvalues and the weights (Corollary \ref{appliGW}). It is natural to state the following conjecture.  

\begin{conj}
\label{GWconjSR}
Let $\mu \in \mathcal M_1(T)$ with Verblunsky coefficients $(\alpha_k)_{k \geq 0}\in {\mathbb{D}}^{\mathbb N}$. Then for any $\g <-1$ and $\mu \in \Sr_1^{\mathbb T}(\pi-\theta_\g,\pi+\theta_\g)$,
\begin{align}
\label{conjGWg+}
{\mathcal K}(\GW_{\g}| \mu) + \sum_{i=1}^{N^-} \mathcal F_{\g}^+(\lambda_i^+) + \sum_{i=1}^{N^+} \mathcal F_{\g}^- (\lambda_i^-) = H(\g) -  
\g \Re \left(\alpha_0 - \sum_{k=1}^\infty \alpha_k\bar\alpha_{k-1}\right) - \sum_{k=0}^\infty \log (1 - |\alpha_k|^2)\,,
\end{align}
where $H$ is defined in (\ref{defH}). 
If $\mu \notin \Sr_1^{\mathbb T}(\pi-\theta_\g,\pi+\theta_\g)$, the right hand side equals $+\infty$.
\end{conj}

This statement would be a direct consequence of Corollary \ref{appliGW} and of Conjecture \ref{GWconjLDP}, as soon as the latter is true. 

\subsection{Higher-order  sum rules}
Besides the Simon sum rule, extensions to higher-order have been tried, either in the research of sum rules or in the research of gems, when the 
 reference measure is 
\begin{equation}
\label{defsigma}
 d\sigma (z) =q(e^{i\theta})  \frac{d\theta}{2\pi}, \qquad \qquad  
\end{equation}
with full support.
For example, Denisov and Kupin  \cite{denisovkupin} considered a reference measure 
$\sigma$ on $\mathbb T$ defined by
\begin{equation}
 q(z) =  \frac{1}{K_q}\prod_1^r |z-\zeta_k|^{2\kappa_k}\,
\end{equation}
where $\zeta_k \in \mathbb T$ and $\kappa_k \in \mathbb N$ (Theorem 2.3 therein). 
The authors used the CMV representation (see \eqref{favardfini}) of operators and give a sum rule. 
Let us call $\mathcal C_0$ the CMV matrix corresponding to $\frac{d\theta}{2\pi}$ (with $\alpha_0 = 1, \alpha_k =0 , k > 0$) 
and $\mathcal C_\sigma$ the operator corresponding to $\sigma$.
If $\mu$ is such that 
\begin{equation}
\label{rank}
\hbox{rank} \ (\mathcal C_\mu - \mathcal C_0) < \infty\,,\end{equation} which means that $\alpha_k  =0$ for $k$ large enough , then 
the sum rule is of the form
\begin{equation}
\label{DenisovKupin}
\mathcal K(\sigma\!\ |\!\ \mu) 
= \mathcal K(\sigma | \operatorname{UNIF}) 
 - \sum_{k \geq 0} \log (1- |\alpha_k|^2) - \Re\!\ \tr \left(Q(\mathcal C_\mu) - Q(\mathcal C_0)\right)
\end{equation}
with $Q$ a polynomial. 
The simplest example is the strongly coupled Gross-Witten model (see Proposition \ref{5.3}), where \[R=1\ , \ q(e^{i\theta}) = 1 - \cos \theta\ , \ a_1= -1/2\  , \ a_k = 0 
\ (k >  2)\ , \ Q(z) =  - z\,,\]  and 
\[- \Re\!\ \tr \left(Q(\mathcal C_\mu) - Q(\mathcal C_0)\right) = (\Re\!\  \tr \!\ \mathcal C_\mu)  - 1\,.\]
This is  in accordance with the Simon sum rule. 

After claiming the sum rule under the finite rank assumption, the authors gave a short proof of an extension to operators satisfying a condition rather hard to check.
For another expression of the sum rule in terms of Verblunsky coefficients with the only  assumption that $q$ in  (\ref{defsigma}) is the square of a polynomial, see \cite{golinskii2007coefficients} Theorem 3.3. For instance, if 
\begin{equation}
\label{antipodal} q(e^{i \theta}) = \frac{1}{K_q} (1 - \cos \theta)^m\,,
\end{equation}
then there exists a function $g$ such that
\[\mathcal K(\sigma | \mu) = \mathcal K(\sigma |\operatorname{UNIF}) 
+ \sum_{k=0}^\infty g( \alpha_{k-m} , \dots, \alpha_k)\]
(see  formulas (2.1), (2.2) in \cite{lukic2013higher} where there is an application for gems).
Another example corresponding to
\begin{equation}
\label{nonantipodal}q(e^{i \theta}) = \frac{1}{K_q} \left(1 - \cos(\theta-\theta_1) \right)  \left(1 - \cos(\theta-\theta_2) \right)\end{equation}
is treated in 
 \cite{simon2005higher}. We also refer to the works \cite{eichinger2014killip} and \cite{Nazarov} for more interesting extensions. 

\section{Matrix extensions}
\label{smamat}
In this section we show how several results can be extended to the case of operator valued measures. 
Since the proofs are mostly identical to the scalar case or can be found in the companion paper \cite{GaNaRomat}, we omit most of them.
In what follows,  $ p$ is a fixed integer $( > 1)$ and the $p\times p$ identity matrix is denoted
by  $\mathbf{1} = \mathrm{I}_p$.

\subsection{Matrix spectral measures}

A matrix measure $\Sigma=(\Sigma_{i,j})$ of size $p\times p$   on $\mathbb{T}$ is a matrix of signed complex measures, such that for any Borel set 
$A\subset \mathbb{T}$, $\Sigma(A)= (\Sigma_{i,j}(A))\in \mathcal H_p$ is (Hermitian and) non-negative definite. A matrix measure on $\mathbb{T}$ is normalized, if $\Sigma(\mathbb{T})=\mathbf 1$. We denote by $\mathcal{M}_{p,1}(T)$ the set of normalized $p\times p$ matrix measures with support in $T\subset\mathbb{T}$. Given a unitary operator $U$ and a collection of vectors $e_1,\dots ,e_p$ cyclic for $U$, one can define the spectral matrix measure $\Sigma$ of $(U,e_1,\dots e_p)$ similarly to (\ref{defspectralmeasure}) by the relation
\begin{align} \label{def:spectralmatrixm}
\langle e_i,U^ke_j\rangle = \int_{\mathbb{T}} z^k \, 
d\Sigma_{i,j}
,\qquad k\in \mathbb{Z},\; i,j\in\{1,\cdots,p\}.
\end{align}
In fact, if $U\in \mathbb{U}(n)$ with eigenvalues $\lambda_1=e^{\ii \theta_1},\dots ,\lambda_n=e^{\ii \theta_n}$ and $\psi_1,\dots ,\psi_n$ is a corresponding system of orthonormal eigenvectors, the spectral matrix measure is given by
\begin{align} \label{matrixmeasure}
\Sigma^{(n)} = \sum_{k=1}^n W_k \delta_{\lambda_k} ,
\end{align}
where $W_k = \psi_k ^{(p)}(\psi_k^{(p)})^\dagger$. Here, $\psi_k^{(p)}$ is the vector of $\mathbb{C}^p$ consisting of the $p$ first coordinates of $\psi_k$.   
Let $\Sigma \in \mathcal{M}_{p,1}(\T)$ be a \emph{quasi scalar} measure, which means that   $\Sigma = \mathbf{1} \sigma$ where $\sigma \in \mathcal{M}_{1}(\T)$ is a scalar measure. Further,  let $\Pi$ be a normalized matrix measure with Lebesgue decomposition
\[\Pi(dz) = h(z) \sigma(dz) +  \Pi^s(dz)\,.\] 
Then, we define
\begin{equation}
\label{kullback}
\mathcal K(\Sigma | \Pi) := - \int_{\mathbb{T}} \log\det h(z)\ \sigma (dz)\,.
\end{equation}
Note that if we define a density matrix componentwise, i.e.
\begin{align*}
\left(\frac{d\Sigma}{d\Pi}\right)_{i,j}  =\frac{d\sigma}{d\Pi_ {i,j}} \ \ , i,j = 1, \dots, p\,,
\end{align*}
 then it is possible to 
rewrite the above quantity in terms of
 the Kullback-Leibler information (or relative entropy) 
\[\mathcal K(\Sigma | \Pi) = \int_{\mathbb{T}} \log \det \frac {d\Sigma}{d\Pi} (z) d\sigma(z)\,,\]
if the density $\frac{d\Sigma}{d\Pi}$ exists and infinity otherwise (see \cite{mandrekar1} or \cite{robertson}). 

As in the scalar case, we can define a matrix version of the Verblunsky coefficients. Now, the construction uses matrix orthogonal polynomials on the unit circle (MOPUC).
\subsection{MOPUC}
\label{MOPUC}

We follow the notation of \cite{FGARop} and \cite{damanik2008analytic}.  A $p\times p$ matrix polynomial $\mathbf{F}$ is a polynomial with coefficents in $\mathbb{C}^{p\times p}$. Given a measure $\Sigma \in \mathcal{M}_{p,1}(\T)$, we define two inner products on the space of $p\times p$ matrix polynomials by setting
\begin{eqnarray*}\lal \mathbf{F}, \mathbf{G}\rar_R &=& \int_\mathbb T \mathbf{F}(z)^\dagger d\Sigma(z)\mathbf{G}(z) \in \mathbb C^{p\times p}\, ,\\
\lal \mathbf{F}, \mathbf{G}\rar_L &=& \int_\mathbb T \mathbf{G}(z) d\Sigma(z)\mathbf{F}(z)^\dagger \in \mathbb C^{p\times p}\, .
\end{eqnarray*}
A sequence of matrix polynomials $(\fbold_j)$ is called right-orthonormal if, and only if, 
\[\lal\fbold_i, \fbold_j\rar_R= \delta_{ij}\mathbf 1 \, .\]
%
As in  the scalar case, we can construct orthonormal polynomials satisfying a recursion and the matrices appearing in this recursion are the so-called matrix Verblunsky coefficients (see \cite{damanik2008analytic} and the historical introduction therein). For the sake of completeness, we give some more details. 
First, assume that the support of $\Sigma$ is infinite. We define the right monic matrix orthogonal polynomials $\Fbold_n^R$ by applying the block Gram-Schmidt algorithm  to $\{\mathbf 1, z\mathbf 1, z^2 \mathbf 1, \dots\}$.  In other words, $\Fbold_k^R$ is the unique matrix polynomial 
$\Fbold_k^R(z) = z^k \mathbf 1 +$ lower order terms, such that $\lal z^j\mathbf 1, \Fbold_k^R\rar_R =0$ for $j=0, \dots, k-1$. 
The normalized orthogonal polynomials are defined by
\[\fbold_0 = \mathbf 1\ \  ,\ \  \fbold_k^R = \Fbold_k^R\kappa_k^R.\]
Here the sequence of $p\times p$ matrices $(\kappa_k^R)$ 
satisfies, for all $k$, the condition $\left(\kappa_k^R\right)^{-1}\kappa_{k+1}^R>0_p$
and is such that the sequence $(\fbold_k^R)$ is orthonormal.
We define the sequence of left-orthonormal polynomials $(\fbold_k^L)$ in the same way except that the above condition is replaced by  $\kappa_{k+1}^L \left(\kappa_k^L\right)^{-1}> 0$. The matrix Szeg\H{o} recursion is then 
\begin{eqnarray}
\label{SzL}
z\fbold_k^L -\rbold_k^L\fbold_{k+1}^L &=& \abold_k^\dagger(\fbold_k^R)^*\\
\label{SzR}
z\fbold_k^R - \fbold_{k+1}^R\rbold_k^R &=& (\fbold_k^L)^* \abold_k^\dagger\,,
\end{eqnarray}
where for all $k\in\mathbb{N}_0$,  
\begin{itemize}
\item $\abold_k$ belongs to  $\mathbb{B}_p$, the closed unit ball of $\mathbb C^{p\times p}$  defined by
\begin{equation}
\label{defBp}\mathbb{B}_p:=\{M \in \mathbb C^{p\times p}: MM^\dagger \leq \mathbf 1\}\,,\end{equation}
\item $\rbold_k$ is the so-called defect matrix defined by
\ben
\label{defrho}
\rbold^R_k :=  \left(\mathbf 1 - \abold_k\abold_k^\dagger\right)^{1/2}\  , \  \rbold^L_k =    \left(\mathbf 1 - \abold_k^\dagger\abold_k\right)^{1/2}\,,
\een
\item for a matrix polynomial $\mathbf{P}$ with degree $k$, the reversed polynomial $\mathbf{P}^*$ is defined by
\[\mathbf{P}^*(z) := z^k \mathbf{P}(1/\bar z)^\dagger\,.\]
\end{itemize}
Notice that  the construction of the recursion coefficients uses only the matrix moments of the matrix measure.
Verblunsky's theorem (the analogue of Favard's theorem for matrix orthogonal polynomials on the unit circle) establishes a one-to-one correspondance between matrix measures on $\mathbb T$ with infinite support and sequences of elements in the interior of $\mathbb B_{p}$ (Theorem 3.12 in \cite{damanik2008analytic}). 

Now, for a matrix measure having a finite support, the construction of the Verblunsky coefficients is not obvious. In \cite{dewag09} Theorem 2.1, a sufficient condition on the moments for such a construction is provided. It is related to the positivity of a block-Toeplitz matrix, as it is also mentioned  in \cite{simon05} at the top of p. 208.
\subsection{Deformed Verblunsky coefficients}
\label{sec:deformedVer}

This section is devoted to a detailed study of  the deformed Verblunsky coefficients in the matrix setting, consisting in identification of  their different definitions and properties. To make the reading easier, we recall the essential results of the scalar case proved in \cite{BNR}. 
\subsubsection{Scalar case}
Motivated by the study of the (scalar) Hua-Pickrell ensemble, Bourgade et al{.} \cite{BNR} introduced the so-called {\it deformed Verblunsky coefficients}. They could be defined in various ways. 
\paragraph{OPUC recursion and the Schur machinery}
Let us assume that $\mu \in \mathcal M_1(\mathbb T)$ has either a finite support consisting of $n$ points, or infinite support and we will say $n =\infty$ and $k \leq n-1$ will mean $k \geq 0$. Then, starting with the monic orthogonal polynomials $\phi_k$ in $L^2(\mu)$ we define 
for $k \leq n-1$ the functions
\begin{align}
\label{defbscal}b_k(z) &:= \frac{\phi_k(z)}{\phi_k^*(z)}  \,, \\
\gamma_k(z) &:= z - \frac{\phi_{k+1}(z)}{\phi_k(z)}\,.
\end{align}
From the Szeg\H{o} recursion (\ref{recpolycirc}), we have 
\[\gamma_k(z) = \frac{\bar\alpha_k}{b_k(z)}\]
and recursively  
\begin{align}
\phi_k(z) &= \prod_0^{k-1} (z - \gamma_j(z)) \ ,\\
\gamma_k(z) &= \bar\alpha_k \prod_0^{k-1} \frac{1- z\widetilde{\gamma}_j(z)}{z- \gamma_j (z)}\ ,\quad  \hbox{with} \ \ 
\widetilde{\gamma}_j(z) =\overline{\gamma_{j}(\bar z^{-1})}\,.
\end{align}
The  deformed Verblunsky coefficients are  by definition 
\[\gamma_k := \gamma_k(1)\ \,\] 
 and may be computed recursively  as 
\[\gamma_0 = \bar\alpha_0 \ , \ \gamma_k = \bar\alpha_k\prod_{j=0}^{k-1} \frac{1- \bar \gamma_j}{1- \gamma_j}\,.\]
If $n$ is finite, the obvious relation $\phi_n (1) = \prod_0^{n-1} (1 - \gamma_j)$ may be lifted up, when $(U,e)$ is given, as 
 \[\det(\mathrm I - U ) = \phi_n (1) =  \prod_{j=0}^{n-1}(1-\gamma_j)\,.\] 

To explain the connection with Schur parameters, let us recall that the Carath\'eodory function of a measure $\mu \in \mathcal M_1(\mathbb T)$
 is defined as 
\begin{equation}
\label{defCara}
F(z) = \int \frac{e^{i \theta} + z}{e^{\ii \theta} -z}d\mu(e^{\ii \theta})\
\end{equation}
and its Schur function $f : \mathbb D  \rightarrow \mathbb D$ is defined through $F$ by:
\[f(z) = \frac{1}{z}\frac{F(z) - 1}{F(z) +1}\,.\]
The Schur algorithm allows to parametrize the Schur function $f$ by a sequence of so-called Schur parameters. 
For $\alpha \in \mathbb D$, let 
\begin{equation}
T_\alpha : \zeta \mapsto (\zeta -\alpha) (1 - \bar\alpha \zeta)^{-1}\,,
\end{equation}
The reverse mapping is $T_{-\alpha}$. If we define the sequence
\begin{equation}
f_0(z) = f(z) \ , \ f_{j+1}(z) = z^{-1} T_{\alpha_j} (f_j(z))  
\ , \ \alpha_j= f_j(0)\,,
\end{equation}
and we say that $f$ is the Schur function associated with the sequence $(\alpha_0, \alpha_1, \dots)$.

The Geronimus theorem states that these are exactly the Verblunsky coefficients.
From the basic recursions
\begin{eqnarray}
\nonumber
\phi_{n+1}(z)  &=& z \phi_n (z) - \bar\alpha_n \phi_n^*(z)\\ 
\label{Schurbasic}
\phi_{n+1}^* (z) &=& \phi_n^* (z) - \alpha_n z \phi_n(z)\,,
\end{eqnarray}
we deduce that the  sequence of quotients $b_k(z)$ defined in (\ref{defbscal}) 
satisfies the recursion
\begin{equation}
b_k(z) = \frac{zb_{k-1}(z) - \bar\alpha_{k-1}}{1- z\alpha_{k-1}b_{k-1}(z)}\,,
\end{equation}
i.e. 
\[b_{k-1}(z) = z^{-1} T_{-\bar\alpha_{k-1}}(b_k(z))\,.\]
In other words, 
 $b_k$ is the Schur function corresponding to the reversed sequence $(-\bar\alpha_{k-1}, \cdots, - \bar\alpha_0, 1)$ (see \cite{Simon2} Prop. 9.2.3), we say that the sequence $(b_k)$ is the sequence of inverse Schur iterates.
\paragraph{Decomposition by reflections}
When $n$ is finite,
 a geometrical interpretation is possible. It relies on the decomposition of $U$ into a product of complex reflections parametrized by the coefficients $\gamma_k,\ k = 0, \dots, n-1$. 

A $n$-(complex) reflection $r$ is an element of $\mathbb U(n)$ such that $r - \mathrm I_n$ has rank $0$ or $1$. If $e$ and $m\not= e$ are unit vectors of $\mathbb C^n$, there is a unique reflection $r$ such that $r(e) =m$, and it is 
\begin{equation}
\label{refl}
r = \mathrm{I}_n - \frac{1}{1- \langle m,e\rangle}(m-e)(m-e)^\dagger
\end{equation}

If $F := \hbox{span}\!\ \{e,m\}$, then $r$ leaves $F^\perp$ invariant. Now setting
\[\gamma = \langle e, m\rangle\ , \  \rho = \sqrt{1 - |\gamma|^2} \ , \ e^{\ii \varphi} = \frac{1- \gamma}{1- \bar\gamma}\,,\]
then, in the basis  $(e,g)$ of $F$ obtained by the Gram-Schmidt procedure, the restriction of $r$  to $F$ has  the matrix
\[\Xi(\gamma) = \begin{pmatrix}\gamma & \rho e^{\ii\varphi}\\\rho& - \bar\gamma e^{\ii \varphi}
\end{pmatrix}\,.\]
Let $U\in \mathbb U(n)$, let $e$ be a cyclic vector for $U$ and let $(\varepsilon_1, \cdots, \varepsilon_n)$ be the orthonormal basis  obtained from the Gram-Schmidt procedure applied to $(e, Ue, \cdots, U^{n-1}e)$. We define recursively $n$ reflections as follows:  $r_1$ is the reflection mapping $e= \varepsilon_1$ onto $Ue = U\varepsilon_1$ and for $k \geq 2$,  $r_k$ is the reflection mapping $\varepsilon_k$ onto
$r_{k-1}^{-1}r_{k-2}^{-1} \cdots r_1^{-1}U \varepsilon_k$. Then $U= r_1 \cdots r_n$ and
\[\langle \varepsilon_k, r_k \varepsilon_k\rangle = \gamma_k\,.\]
\subsubsection{The matrix case}
\label{2.10}
\paragraph{MOPUC recursion and the Schur machinery}
Let us define, for $k = 0, \dots, n-1$,
\begin{equation}
\label{defbmat}
\bbold_0(z) = \mathbf{1}, \ \ , \ \bbold_k (z) = \fbold_k^L(z) \left(\fbold_k^{R, *}(z)\right)^{-1}\,,
\end{equation}
(notice this is consistent with the definition (\ref{defbscal}) since when $p=1$, $\kappa_n^R$ and $\kappa_n^L$ are scalar and equal). We also set
\begin{equation}\label{phioverphimat}
\bbold_k = \bbold_k(1) \ 
 \ \hbox{and} \ \ \gbold_k = \bbold_k^{-1}\abold_k^\dagger\,.
\end{equation}
These coefficients $\gbold_k$ are called  deformed matrix Verblunsky coefficients. 
 
As in the scalar case, we can make the connection with the inverse Schur iterates. 

The Carath\'eodory function $\bf{F}$ is now matrix-valued, defined again by (\ref{defCara}), and the Schur function is (\cite {damanik2008analytic} Prop. 3.15)
\[\mathbf{f}(z) = z^{-1}(\mathbf{F}(z) - \mathbf{1})(\mathbf{F}(z) + \mathbf{1})^{-1}\,.\]
To define the Schur algorithm, we set for $\abold \in \mathbb B_p$ with $\abold \abold^\dagger< \mathbf{1}$,
\[\mathbf{T}_{\abold} (\boldsymbol{\zeta})= (\rbold^R)^{-1}(\boldsymbol{\zeta} - \abold)(\mathbf{1} - \abold^{\dagger}\boldsymbol{\zeta})^{-1}\rbold^L\]
The reverse mapping is $\mathbf{T}_{-\abold}$, and we notice that
\begin{equation}
\label{magicT}
\left(\mathbf{T}_{\abold}(\zeta)\right)^{-1} = \mathbf{T}_{\alpha^\dagger}(\zeta^{-1})\,.
\end{equation}

\begin{prop}[\cite{damanik2008analytic} Th. 3.19]
\label{6.1}
For the Schur functions $\mathbf{f}_0, \mathbf{f}_1, \dots$ associated with Verblunsky coefficients $\abold_0, \abold_1, \dots$, the following relations hold:
\begin{eqnarray}
\mathbf{f}_{j+1}(z) &=& z^{-1} \mathbf{T}_{\abold_j}(\mathbf{f}_j(z))\\
\abold_j &=& \mathbf{f}_j(0)
\end{eqnarray}
\end{prop}

The connection with \eqref{defbmat} is the following.

\begin{prop}[\cite{damanik2008analytic} Prop. 3.26]
\label{6.2}
For $k \geq 1$, $\mathbf{b}_k(z)$ is the Schur function associated with the reversed sequence 
$(- \abold_{k-1}^\dagger, \dots, - \abold_0^\dagger, \mathbf{1})$.
\end{prop}

\paragraph{Decomposition by reflections}

Let us first fix some more notations. Let $\mathbf{e} = [e_1, \dots, e_p]$ be a $N \times p$ matrix  consisting of $p$ column vectors of dimension $N=np$. If $U \in \mathbb U(N)$, we denote by $U\mathbf{e}$ the $N \times p$ matrix $U\mathbf{e} := [Ue_1, \dots, Ue_p]$ .   
The pseudo-scalar product of $\mathbf{e}$ with  $\mathbf{f} = [f_1 \cdots f_p]$ is a $p\times p$ matrix denoted by $\ll \mathbf{f}, \mathbf{e}\gg$ and  defined by
\[\ll \mathbf{f}, \mathbf{e}\gg_{i,j}\!\ = \langle f_i , e _j\rangle \ \  i,j = 1, \dots, p\,.\]
Assume that $\mathbf{e}$ is cyclic for $U$ (see definition 2.3 in \cite{FGARop}).  
For $N = np$ with $n \geq 2$, 
let $(\boldsymbol{\varepsilon}_1, \dots, \boldsymbol{\varepsilon}_Q)$  be the orthonormal basis obtained from the Gram-Schmidt procedure applied to
$(\mathbf{e}, U\mathbf{e}, \dots, U^{Q-1}\mathbf{e})$.  
Neretin \cite{Ner1} defined a sequence of operations on unitary matrices of decreasing dimensions recalled here.
For $ m< n$ we decompose a unitary matrix $U \in \mathbb U(n)$ into four blocks
\[U = \begin{pmatrix} A& B \\
C&D\end{pmatrix}\] 
with $A$ a $m\times m$ matrix, and then define
\[\Xi_n^m (U) = D + C ({\mathrm I}_m-A)^{-1} B \in \mathbb U(n-m)\,. \]
Actually, ${\mathrm I}_{n-m} - \Xi_n^m (U)= ({\mathrm I}_n -U) \diagup  ({\mathrm I}_m-A)$ where $M \diagup N$ is the Schur complement of $M$ with respect to its upper left block (submatrix) $N$.
 This doubly indexed sequence of transformation enjoys the projective property:
\begin{equation}
\label{projxi}\Xi_{n-m}^r \circ \Xi_n^m = \Xi_n^{r+m}\,,\end{equation}
as soon as $r+m < n$ (see \cite{Ner1}, Proposition 0.1). In the sequel, for $q > p$, we denote by $[M]_p$ the upper left block of the $q \times q$ matrix $M$.

We define the successive iterations 
\begin{equation}
\label{defcbold}
\mathbf{c}_0(U) := [U]_p \ , \ \mathbf{c}_r(U) := [\Xi_N^{rp}(U)]_p  \   \ , \ 1 \leq r \leq n-1\,.
\end{equation}
Then Neretin proved (\cite{Ner1}, Section 1.5)
\begin{equation}
\label{capital}
\det ({\mathrm I}_N - U) = \prod_{r=0}^{n-1}\det \left( \mathbf{1} - \mathbf{c}_r(U) \right)\,.
\end{equation}
These operators $\Xi$ are  
used to define the successive reflections. More precisely, we define  
\[\hat\pi_0 (U) = U \ , \ \hat\pi_k (U) = {\mathrm I}_{kp} \oplus \Xi_N^{kp}(U) \ , \ 1 \leq k \leq n-1 \ , \ \hat\pi_Q (U)= {\mathrm I}_N\,\]
\[R_j(U) = \hat\pi_{j-1}(U) \hat\pi_j(U)^\dagger\ , \  1 \leq j \leq n\,.\]
If $U$ is written using an orthonormal  basis $(\mathbf{e}, \mathbf{e}_2, \dots, \mathbf{e}_n)$, then  $R_1$ maps $\mathbf{e}$ onto $U\mathbf{e}$ and is a reflection since the rank of $R_1 - {\mathrm I}_N$ is the same as the rank of $U - ({\bf 1} \oplus \Xi_N^p)$ which is at most $p$ (see Prop. 2.5 in  \cite{bourgade2012unitary}).

More generally, for $k \geq 2$, $R_k$ is 
a reflection mapping $\boldsymbol{e}_k$ onto $R_{k-1}^\dagger R_{k-2}^\dagger \dots R_1^\dagger U\boldsymbol{e}_k$ and

 \[U= R_1 \dots R_Q\,.\]
In particular, let $\mathcal G$ be the matrix of a unitary operator $U$ written in the basis $(\boldsymbol{\varepsilon}_k)$ obtained by  orthonormalizing the sequence $\mathbf{e}, U\mathbf{e}, \dots, U^{n-1}\mathbf{e}$.  Usually  $\mathcal G$ is called the block GGT matrix : 
\begin{eqnarray*}
\mathcal G  := \mathcal G^R(\abold_0, \abold_1, \dots) = \begin{pmatrix} \abold_0^\dagger &\rbold_0^L\abold_1^\dagger&\rbold_0^L \rbold_1^L \abold_2^\dagger& \rbold_0^L\rbold_1^L \rbold_2^L \abold_3^\dagger & \dots\\
\rbold_0^R&-\abold_0\abold_1^\dagger& -\abold_0\rbold_1^L\abold_2^\dagger&- \abold_0\rbold_1^L \rbold_2^L \abold_3^\dagger& \dots\\
0 & \rbold_1^R&  - \abold_1\abold_2^\dagger&- \abold_1\rbold_2^L\abold_3^\dagger& \dots\\
0&0 &\rbold_2^R & - \abold_2\abold_3^\dagger & \dots\\
\vdots&\vdots&\vdots& \vdots&\ddots
\end{pmatrix}
\end{eqnarray*}
Write $\mathcal R_j$ for $R_j(\mathcal G)$. Then $\mathcal R_1$ maps $\boldsymbol{\varepsilon}_1$ onto $\mathcal G \mathbf{\varepsilon}_1$ and is a reflection since the rank of $\mathcal R_1 - {\mathrm I}_N$ is the same as the rank of $\mathcal G - ({\bf 1} \oplus \Xi_N^p)$ which is at most $ p$ (see Prop. 2.5 in  \cite{bourgade2012unitary}).

More generally, for $k \geq 2$,  $\mathcal R_k$ is 
 a reflection mapping $\boldsymbol{\varepsilon}_k$ onto $\mathcal R_{k-1}^\dagger\mathcal R_{k-2}^\dagger \dots \mathcal R_1^\dagger \mathcal G\boldsymbol{\varepsilon}_k$ and
 \[\mathcal  G= \mathcal R_1 \dots \mathcal R_Q\]
Of course, we have
\begin{equation}
\label{detfond}
\det (\mathrm{I}_N -U) = \det (\mathrm I_N - \mathcal G) = \prod_{r=0}^{n-1}\det (\mathbf{1} - \mathbf{c}_r (\mathcal G)),\end{equation}
and we have the following identification.
\begin{prop}
\label{identif}
Let $(U, \bf{e})$ be given and call $\mathcal G$ the matrix of $U$ written in the basis $(\boldsymbol{\varepsilon}_k)$. Then for $k=1, \dots, n$
\begin{eqnarray}
\label{new}
\mathbf{c}_{k-1} (\mathcal G) =\!\ \ll \boldsymbol{\varepsilon}_k , \mathcal R_k \boldsymbol{\varepsilon}_k \gg\!\ = \gbold_{k-1}\,. 
\end{eqnarray}
\end{prop}

\subsection{LDP for matrix Verblunsky coefficients}

In a previous work (\cite{FGARop}), the first and last author studied the  CUE case. If $N = np$, the  matrix Verblunsky coefficients $\abold_0, \dots, \abold_{n-1}$ are independent, and for $k \leq n-2$, $\abold_k$ has  a density in $\mathbb B_p$  proportional to
\begin{equation}
\label{mV0}
 \det (\mathbf 1 - \abold\abold^\dagger)^{(n-k-2)p}\,,\end{equation}
 (it is a matricial extension of $\eta_{n-k-2}$ defined in (\ref{defeta})). Note that all densities involved in this section are with respect to 
\begin{align*}
dM = \prod_{1\leq k , l\leq p} d(\Re M_{k l}) \!\  \prod_{1\leq k , l\leq n} d(\Im M_{kl}) .
\end{align*}
From this density, we deduced the LDP:

\begin{prop}[\cite{FGARop} Theorem 3.6]
\label{lemmatunit} For $N=np$, let $U\in \mathbb{U}(N)$ be drawn 
from the Haar measure $\mathbb{P}^{(N)}$. Let further $(\abold_k^{(n)})_{0\leq k\leq n-1}$ be the matrix Verblunsky coefficients of the spectral matrix measure of $(U,e_1,\dots e_p)$. Then, for any fixed $k\geq 1$,  $(\abold_0^{(n)}, \abold_1^{(n)}, \cdots, \abold_k^{(n)})_{n \geq k}$ satisfies 
the LDP in $(\mathbb B_p)^k$ with speed $N$ and good rate function
\[I_k(\abold_0, \dots, \abold_k) = \sum_{j=0}^k  -\log\det (\mathbf 1 - \abold_j \abold_j^\dagger) \,.\]
\end{prop}

To study the Hua-Pickrell case, we will use  
the deformed matrix Verblunsky coefficients defined in (\ref{phioverphimat}). 
%
Their distribution is given in the following proposition, whose proof is postponed to Section \ref{sec:alphajlaw}.


\begin{thm}
\label{alphajlaw}
Let $N = np$ with  $n > 2$ and $U \in \mathbb U(N)$ be drawn from the Hua-Pickrell 
 distribution $\mathbb H \mathbb P^{(N)}_\delta$
. Let $(\gbold_k^{(n)})_{0\leq k\leq n-1}$ be the deformed matrix Verblunsky coefficients of the spectral matrix measure of $(U,e_1,\dots e_p)$. Then, $\gbold_1\sn, \dots ,\gbold_{n-1}\sn$ are independent. 
Moreover, for $k \leq n-2$, $\gbold_k\sn$ has in $\mathbb B_p$ the density 
\begin{equation}
\label{densalphajmat}
K_{n,k}^{(\delta)} \det \left(\mathbf 1  - \gbold\right)^{\bar\delta}\det \left(\mathbf 1  - \gbold^\dagger\right)^\delta \det (\mathbf 1 - \gbold\gbold^\dagger)^{(n-k-2)p}
\end{equation}
where
\begin{equation}
\label{cstner}
K_{n,k}^{(\delta)} = \pi^{- p^2}\prod_{j=1}^p\frac{\Gamma(N- (k+1)p +j + \delta) \Gamma(N- (k+1)p +j + \bar \delta)} {\Gamma(N-(k+2)p +j) \Gamma(N - (k+1)p +j+ \delta+\bar\delta)}
\end{equation}
and $\gbold_{n-1}$  follows the Hua-Pickrell distribution on $\mathbb U (p)$ with parameter $\delta$.
\end{thm}

If $\delta=N\d$, we get the following LDP for the deformed coefficients. We remark that if $\d=0$, the rate function is that of  Proposition \ref{lemmatunit}.  Indeed, the matrices $\bbold_k$ are unitary by Theorem 3.9 in \cite{damanik2008analytic}. 

\begin{prop}
\label{lemmat}
Let $N = np$ with  $n > 2$ and $U\in \mathbb U(N)$ be drawn  from the Hua-Pickrell probability distribution  $\mathbb H\mathbb P_{N\d}^{(N)}$ ($\d\geq 0$). Let $(\gbold_k^{(n)})_{0\leq k\leq n-1}$ be the deformed matrix Verblunsky coefficients of the spectral matrix measure of $(U,e_1,\dots e_p)$. 
Then, for any fixed $k$, $(\gbold_0\sn, \gbold_1\sn, \cdots, \gbold_k\sn)_{n \geq k}$ satisfies the LDP in $(\mathbb B_p)^k$ with speed $N$ and good rate function
\[I_k(\gbold_0, \dots, \gbold_k) = \sum_{j=0}^k H_{\d, p} (\gbold_j)\,,\]
with
\begin{equation}
\label{defHdp}
H_{\d, p}(\gbold) = -\log\det (\mathbf 1 - \gbold \gbold^\dagger) - \d \log \det\left((\mathbf 1 - \gbold)(\mathbf 1 - \gbold)^\dagger\right) + p H_\d (0)\end{equation}
where $H_\d(0)$ is defined in (\ref{H1}).
\end{prop}
Similarly to the scalar case, the function  $H_{\d,p}$ is nonnegative and vanishes uniquely at $\gbold = \gamma_\d \cdot \mathbf 1$.  
\subsection{LDP for matrix spectral measures}

Our next LDP holds for matrix spectral measures of $(U,e_1,\dots ,e_p)$, when $U$ is drawn with the general measure $\mathbb P^{(N)}_\V $ on $\mathbb{U}(N)$ as defined in \eqref{PnV}. In this case, the eigenvector matrix is again Haar distributed, so that the weights $(W_1,\dots ,W_N)$ are independent of the eigenvalues. Moreover, they follow a distribution that is a  matrix analogue of the Dirichlet law. For the precise statement, we refer to Proposition 3.1 in \cite{FGARop}.
Let us introduce the matrix analogue of the set $\mathcal{S}^{\mathbb T}_1 = \mathcal{S}^{\mathbb T}_1(\am,\ap)$. For $[\alpha^- , \alpha^+]$ an interval included in $(0, 2\pi)$, let 
$I=\widehat{[\alpha^- , \alpha^+]}$
and let $\mathcal{S}^{\mathbb T}_{p,1} = \mathcal{S}_{p,1}^{\T}(\am,\ap)$ be the set of all normalized measures 
$\Sigma \in\mathcal{M}_{p,1}(\mathbb{T})$ with 
\begin{itemize}
\item[(i)] $\operatorname{supp}(\Sigma) = J \cup \{e^{\ii \theta_i^-}\}_{i=1}^{N^-} \cup \{e^{\ii \theta_i^+}\}_{i=1}^{N^+}$, where $J\subset I$, $N^-,N^+\in\N\cup\{\infty\}$ and $\theta_i^\pm \in [0,2\pi)$. Furthermore,
\begin{align*}
0\leq\theta_1^-<\theta_2^-<\dots <\am \quad \text{and} \quad \theta_1^+>\theta_2^+>\dots >\ap .
\end{align*}
\item[(ii)] If $N^-$ (resp. $N^+$) is infinite, then $\theta_i^-$ converges towards $\am$ (resp. $\theta_i^+$ converges to $\ap$).
\end{itemize}
We can write such a measure $\Sigma$ as
\begin{align}
\label{muinS1}
\Sigma = \Sigma_{|I} +  \sum_{i=1}^{N^+} \Gamma_i^+ \delta_{\lambda_i^+} + \sum_{i=1}^{N^-} \Gamma_i^- \delta_{\lambda_i^-},
\end{align}
for some nonnegative Hermitian matrices $\Gamma_1^+,\cdots, \Gamma_{N^+}^+,\Gamma_1^-,\cdots, \Gamma_{N^-}^-$ and $\lambda_i^\pm =e^{\ii \theta^\pm_i}$. As before, $\Sr_{p,1}^{\mathbb T}(0,2\pi)$ is the extended notation for the case of matrix measures supported by $\mathbb T$. The proof of the following result is omitted. The steps to extend the scalar case to the matrix case  in Theorem \ref{MAINT} are exactly the same as in \cite{GaNaRomat}. Therein, the  LDP for matrix measures on the real line is established.

\begin{thm}\label{MAINTmatrix}
Assume that $U$ is distributed according to $\mathbb{P}_{\V}^{(N)}$, $N=np$, and that the potential $\V$ satisfies assumptions (T1), (T2), (T3). Then the sequence of matrix spectral measures $(\Sigma^{(N)}_p)_n$ of $(U,e_1,\dots e_p)$ satisfies the LDP in $\mathcal M_{p,1}(\mathbb{T})$ equipped with the weak topology, with speed $N$ and rate function
\begin{equation}
\mathcal I_\mathcal V^p (\Sigma) = \mathcal K (\mathbf 1\cdot \mu_\mathcal V \!\ | \!\ \Sigma ) + \sum_{i=1}^{N^+} \mathcal F_{\V}^+(\lambda_i^+) + \sum_{i=1}^{N^-} \mathcal F_{\V}^-(\lambda_i^-)\,,
\end{equation} 
if $\Sigma \in \Sr_{p,1}^{\mathbb T}(\alpha^-,\alpha^+)$ and $\mathcal I_\mathcal V^p (\Sigma)=+\infty$ otherwise. Here, $\mu_\V$ is the scalar measure as in assumption (T2). 
\end{thm}

\subsection{Sum rules}

The matrix version of  Szeg\H{o}'s formula was established in \cite{delsarte} (see more recently \cite{derevy}). For a probabilistic point of view, see \cite{FGARop}.

\begin{thm} \label{classicalsumrule}
Let $\Sigma \in \mathcal M_{p,1} (\mathbb T)$ with infinite support and let $(\abold_k)_{k \geq 0}  \in (\mathbb B_p)^{\mathbb N_0}$ be the sequence of its Verblunsky coefficients. Then 
\ben
 \mathcal K(\mathbf 1\cdot \operatorname{UNIF}| \Sigma)
= \sum_{k=0}^\infty - \log \det (\mathbf 1 - \abold_k \abold_k^\dagger)\,. 
\een
\end{thm}
Our next result is a matrix version of Theorem \ref{sumruleHP}. It is a combination of Proposition \ref{lemmat} and Theorem \ref{MAINTmatrix} ($\V$ is here the potential of the Hua-Pickrell ensemble). The proof follows that of the scalar case given in Section \ref{sec:proofs}: Proposition \ref{lemmat} yields, by the projective method, a complementary LDP for a measure distributed according to the Hua-Pickrell ensemble. Then, the statement follows from the uniqueness of a rate function. Note that for $\d=0$, it reduces to Theorem \ref{classicalsumrule}.

\begin{thm} \label{sumruleHPmatrix}
Let $\Sigma \in \mathcal M_{p,1}(\mathbb T)$ with infinite support and let 
$(\gbold_k)_{k \geq 0} \in (\mathbb B_p)^{\mathbb N_0}$ be the sequence of its deformed matrix Verblunsky coefficients. Then for any $\d\geq 0$, $\sum_{k=0}^\infty H_{\d,p}(\gbold_k)=\infty$ if $\Sigma \notin \Sr_{p,1}^{\mathbb T}(\theta_\d,2\pi-\theta_\d)$. For $\Sigma \in \Sr_{p,1}^{\mathbb T}(\theta_\d,2\pi-\theta_\d)$,
\ben
 \mathcal K(\mathbf 1\cdot \HP_\d | \Sigma) + \sum_{i=1}^{N^+} \mathcal F_{HP}^+(\lambda_i^+) + \sum_{i=1}^{N^-} \mathcal F_{HP}^-(\lambda_i^-) = \sum_{k=0}^\infty H_{\d ,p}(\gbold_k)\,,
\een
where both sides may be infinite simultaneously, and $H_{\d, p}$ is defined in \eqref{defHdp}.
\end{thm}

\begin{cor}
\label{gemHPmatrix}
Let $\Sigma \in \mathcal M_{p,1}(\mathbb T)$ with infinite support and deformed Verblunsky coefficients $(\gbold_k)_{k \geq 0}\in (\mathbb B_p)^{\mathbb N_0}$. Then
\begin{align*}
\sum_{k=1}^\infty \left|\left|\gbold_k - \gamma_\d \mathbf 1 
 \right|\right|^2 < \infty
\end{align*}
if and only if
\begin{enumerate}
\item 
  $\Sigma \in \Sr_{p,1}^{\mathbb T}(\theta_\d,2\pi-\theta_\d)$
\item $\sum_{i=1}^{N^+} (\theta_i^+ - 2\pi+\theta_\d)^{3/2} + \sum_{i=1}^{N^-} (\theta_\d - \theta_i^- )^{3/2}  < \infty$ and if $N^->0$, then $\theta_1^- > 0$. 
\item If $d\Sigma(\theta) = F(\theta)\frac{d\theta}{2\pi}+d\Sigma_s(\theta)$ is the decomposition of $\Sigma$ with respect to the Lebesgue measure, then
\begin{align*}
\int_{\theta_\d}^{2\pi-\theta_\d}  \frac{\sqrt{\sin^2(\tfrac{\theta}{2}) -
\sin^2(\tfrac{\theta_{\d}}{2})}}{2\pi \!\ \sin(\tfrac{\theta}{2})} \log\det(F(\theta)) d\theta >-\infty .
\end{align*}
\end{enumerate}
\end{cor}

The proof of this corollary is very similar to the proof of Corollary \ref{gemHP} (scalar case) and will be omitted.

\medskip

For the Gross-Witten ensemble, it seems difficult (at least at a first attempt) to adapt Simon's proof to the matrix setup. Nevertheless, the density of $\mathbb G\mathbb W^{(N)}_\g$ with respect to $\mathbb P^{(N)}$ 
 is proportional to 
\begin{align*}
\exp \left( N\g \ \Re\!\ \tr \!\ U \right)
\end{align*}
and $\tr\!\ U = \tr\!\ \mathcal G$.
Further, $\tr\!\ \mathcal G$ can be computed in matrix terms, taking into account the GGT  form of $U$:
\[\tr\!\ \mathcal G = \tr\!\ T_n(\abold_0,\dots ,\abold_{n-1})\]
where 
\begin{align*}
 T_n(\abold_0,\dots ,\abold_{n-1}) = \abold_0^\dagger - \abold_0\abold_1^\dagger -\ \dots\ 
 \begin{cases}
 - \abold_{2r}^\dagger\abold_{2r-1}  & \hbox{ if } n=2r+1\,, \\
 -\abold_{2r}\abold_{2r+1}^\dagger & \hbox{ if } n=2r+2 ,\ r\geq 0\,,
 \end{cases}
\end{align*}
and
\begin{align} \label{Talpha}
T(\abold_0,\abold_1,\dots) = \abold_0^\dagger - \abold_0\abold_1^\dagger -\sum_{k=1}^\infty (\abold_{2k}\abold_{2k+1}^\dagger + \abold_{2k}^\dagger\abold_{2k-1}).
\end{align}
We may formulate the matrix version of Corollary \ref{coreasy} and of the Conjecture \ref{GWconjSR}.

\begin{conj}
Let $\Sigma \in \mathcal M_{p,1} (\mathbb T)$ with infinite support and let $(\abold_k)_{k \geq 0} \in (\mathbb B_p)^{\mathbb N_0}$ be the sequence of its matrix Verblunsky coefficients.
\begin{enumerate}
\item
If $|\g| \leq 1$, then
\begin{equation}
\label{summatgwg}
\mathcal K(\mathbf 1 \cdot \GW_{-\g} | \Sigma) = \sum_{k=0}^\infty -\log \det (\mathbf 1 - \abold_k \abold_k^\dagger) +\g \Re\!\  \tr\!\  T (\abold_0,\dots) + H_p(\g)\,,
\end{equation}
where $T (\abold_0,\dots)$ is given by \eqref{Talpha}
and $H_p (\g)$ is some constant.
\item If $|\g| > 1$, then  a similar identity holds, with an additional term on the left hand side which is
\[\sum_{i=1}^{N^-}\mathcal F^+_{-\g}(\lambda_i^-) + \sum_{n=1}^{N^+} \mathcal F^-_{-\g}(\lambda_i^+)\,. \]
\end{enumerate}
\end{conj}
\section{Proofs} \label{sec:proofs}
\subsection{Proofs of Section \ref{sLDDPP}}
\subsubsection{Proof of Theorem \ref{MAINT}}

We can follow verbatim the proof of the corresponding theorem in the real case. The main idea is to apply the projective method (the Dawson-G\"artner Theorem, see \cite{demboz98}) to a non-normalized version of the spectral measure. In a first step, we consider instead of $\mu\sn$ the measure
\begin{align*}
\tilde{\mu}\sn = \sum_{i=1}^n \gamma_i \delta_{\lambda_i}.
\end{align*}
Here, $\gamma_1,\dots ,\gamma_n$ are i.i.d. random variables with distribution $\operatorname{Gamma}(1,n^{-1})$ (their mean is $n^{-1}$). The self-normalized vector built with this sample has a uniform distribution on the simplex. So that,  $\tilde{\mu}\sn/\tilde{\mu}\sn(\T)$ recovers the original distribution of $\mu\sn$. Further, we consider the measure 
\begin{align} \label{projectedmeas}
\pi_j(\tilde{\mu}\sn) = \tilde{\mu}\sn_{|I} +  \sum_{i=1}^{N^+\wedge j} \gamma_i^+ \delta_{\lambda_i^+} + \sum_{i=1}^{N^-\wedge j} \gamma_i^- \delta_{\lambda_i^-}
\end{align}
using the representation as in \eqref{muinS0R}. Note that this projection is not continuous in the weak topology, and in \cite{GaNaRo} we introduce a new topology generated by $\tilde{\mu}\sn_{|I}$ and the vector of outliers. On the set of normalized measures, this topology is stronger than the weak topology and we can claim the LDP in the latter topology. Ultimately, this also explains why our arbitrary distinction between $\lambda_i^+$ and $\lambda_i^-$ creates no problems: the transition of an eigenvalue from $e^{\ii \theta_1^-}$ to $e^{\ii \theta_1^+}$ is continuous in the weak topology, but not in our new one.

A crucial ingredient in the LDP for $\pi_j(\tilde{\mu}\sn)$ is the 
LDP for a finite collection of extreme eigenvalues. Of course, for the first statement of our theorem, this can be omitted. 
For $A$ a subset of $\mathbb{R}$, let $A^{\uparrow j}$ (resp. $A^{\downarrow j}$) denotes the subset of 
$A^j=A\times\cdots\times A$ consisting in all non-decreasing sequences (resp. non-increasing sequences) of $A$.

\begin{prop}
\label{LDPjextremeU}
Let $j$ and $\ell$ be  fixed integers. Assume that $\V$ satisfies (T1), (T2) and the control condition (T3).
If $0 <\alpha^+$ and $\alpha^-<2\pi$, then
the law of $(\theta_1^+,\dots,\theta_j^+,\theta_1^-,\dots, \theta_\ell^-)$ under $\mathbb{P}_\V\sn$ satisfies the LDP in $\R^{j+ \ell}$ with speed $n$ and rate function
\begin{align*}
\mathcal{I}_{\theta^\pm}(\theta^+, \theta^-) = 
\sum_{i=1}^j 
\Fr_{\V}^+(e^{\ii \theta_i^+})  + \sum_{i=1}^\ell 
\Fr_{\V}^-(e^{\ii \theta^-_i})
\end{align*} 
if 
$\theta^+=(\theta_1^+, \dots , \theta_j^+)\in [\alpha^+, 2\pi]^{\downarrow j}$ and $\theta^-= (\theta_1^-, \dots , \theta_\ell^-)\in [0, \alpha^-]^{\uparrow \ell}$ and $\mathcal{I}_{\theta^\pm}(\theta^+, \theta^-) = \infty$ otherwise.
\end{prop}

\proof We first mention the main points in the proof of the large deviation upper bound. 
Let us stress that exponential tightness is inherent on the circle. The proof follows the same lines as in \cite{GaNaRo} and
makes use of the following lemmas.

\begin{lem}
\label{LDPtildeU}
Let $\V$ be  a continuous potential on $\mathbb T\setminus \{1\}$ satisfying (T1) 
 and let $r$ be a fixed integer.  If  $\mathbb P^{(n)}_{\V_n}$ is the probability measure associated to the potential
$\V_n= \frac{n+r}{n} \V$, 
then the law of $\muun$ under $\mathbb P^{(n)}_{\V_n}$ satisfies the LDP with speed $n^2$ and good rate function 
\begin{equation}  
\label{entropyVU}
\mu \mapsto \mathcal E(\mu) - \inf_\nu\mathcal E(\nu)
\end{equation}
where $\mathcal E$  is defined in (\ref{ratemuu}).
\end{lem}

\begin{lem}
\label{teknikU}
If the potential $\V$ is  continuous on $\mathbb T\setminus \{1\}$ and satisfies (T1),
 we have for every $p \geq 1$
\begin{equation}
\lim_{n\to \infty} \frac{1}{n} \log \frac{Z_\V^{(n)}}{Z_{\frac{n}{n-p}\V}^{(n-p)}} = -\inf_{z_1, \dots ,z_k} \sum_{k=1}^p  \mathcal J_\V (z_k) = -p \inf_z \mathcal J_\V(z)\,.
\end{equation}
\end{lem}

\medskip

For the proof of the large deviation lower bound, we may make the same remark as above. We do not need to show exponential tightness anymore. Besides we need the fact that under $\mathbb P^{(n)}_{\frac{n+r}{n}\V}$, the extremal eigenvalues converges to the endpoints of the support of $\mu_V$ if its support is a proper arc. It was a separate lemma in \cite{GaNaRo}, but it is a direct consequence of the upper bound and assumption (T3).

\begin{lem}\label{cvprobaU}
Under Assumptions (T1) and (T3), the distance of $\theta_i^+$ and $\theta_i^-$ to $\{\alpha^-,\alpha^+\}$ converges in probability to 0 for all $i\geq 1$.
\end{lem}

\proof We may use the large deviation upper bound of Proposition \ref{LDPjextremeU}. The upper bound involves the rate function $\mathcal J_\V - \inf_x \mathcal J_\V(x)$. This rate function may vanish somewhere on the support of $\mu_\V$. But, it does not vanish outside of this
support (by assumption (T3)). It follows that the probability that the distance to $\{\alpha^-,\alpha^+\}$ is greater than $\varepsilon$ is exponentially small.\hfill $\Box$

The next step in the proof of Theorem \ref{MAINT} is a joint LDP for the measure $\tilde\mu_{|I(j)}^{(n)}$ restricted to $I(j) = I\setminus \{\lambda^+,\lambda^-\}$ and the extremal eigenvalues. The crucial ingredients are the independence of the eigenvalues and the weights and the LDP for $\mu_\u\sn$ at the faster speed $n^2$. The following result is a straightforward counterpart of Theorem 4.2 of \cite{GaNaRo}.

\begin{prop} \label{ldpjoint} \ 
\begin{enumerate}
\item 
Assume that the potential satisfies (T1) and that the support of $\mu_{\mathcal V}$ is $\mathbb T$. Then the sequence of measures $\tilde{\mu}\sn$ satisfies the LDP with speed $n$ and good rate function
\begin{align*}
\mathcal I(\mu) = \mathcal{K}(\mu_\V\!\ |\!\ \mu) + \mu(\T)-1\,.
\end{align*}
\item
Assume that the potential $\V$ satisfies the assumptions (T1), (T2) and (T3). Then the sequence 
\begin{align*}
( \tilde{\mu}\sn_{|I(j)} , \theta^+, \theta^- )
\end{align*}
with $\theta^\pm= (\theta_1^\pm,\dots ,\theta_j^\pm)$ satisfies the LDP with speed $n$ and good rate function
\begin{align*}
\mathcal{I}(\mu,x^+,x^-) = \mathcal{K}(\mu_\V\!\ |\!\ \mu) + \mu(I)-1
+ \mathcal{I}_{\theta^\pm}(\theta^+, \theta^-) .
\end{align*}
\end{enumerate}
\end{prop}

The weights $(\gamma_1^+,\dots ,\gamma_j^+,\gamma_1^-,\dots ,\gamma_j^-)$ associated with the outlying eigenvalues satisfy the LDP in $\R^{2j}$ with speed $n$ and good rate function 
\begin{align*}
\mathcal{I}_\gamma (y)= 
\begin{cases}
{\displaystyle \sum_{i=1}^{2j} y_i} & \mbox{ if all $y_i$ are nonnegative,}\\
\infty & \mbox{ otherwise.}
\end{cases}
\end{align*}
Using again the independence, we obtain from Proposition \ref{ldpjoint} the joint LDP for 
\begin{align*}
\big( \tilde{\mu}\sn_{|I(j)} , \theta^+, \theta^- , \gamma^+,\gamma^-\big)
\end{align*}
(omitting the outlying eigenvalues if the support of $\mu_\V$ is $\T$). The rate function is the sum of the rate function of Proposition \ref{ldpjoint} and $\mathcal{I}_\gamma$. From this collection, we may now conclude the LDP for the projected measure $\pi_j(\tilde{\mu}\sn)$, by mapping continuously
\begin{align*}
\big( \tilde{\mu}\sn_{|I(j)} , \theta^+, \theta^- , \gamma^+,\gamma^-\big) \longmapsto
\tilde{\mu}\sn_{|I(j)} + \sum_{i=1}^j \left( \gamma_i^+ \delta_{e^{\ii \theta_i^+}} + \gamma_i^- \delta_{e^{\ii \theta_i^-}}\right). 
\end{align*}
It yields by the contraction principle the LDP for $\pi_j(\tilde{\mu}\sn)$ with good rate function 
\begin{align*}
\mathcal{I}_j (\mu) = \mathcal{K}(\mu_\V\!\ |\!\ \mu) + \mu(\T)-1 + \sum_{i=1}^{N^+\wedge j} \Fr_{\V}^+(e^{\ii \theta_i^+}) + \sum_{i=1}^{N^-\wedge j} \Fr_{\V}^+(e^{\ii \theta_i^-}) .
\end{align*}
Finally, the LDP for $\mu^{(n)}$ follows by taking the projective limit and normalizing. The arguments are as in Section 4.4 of \cite{GaNaRo}. This concludes the proof of Theorem \ref{MAINT}.

\subsubsection{Proof of Theorem \ref{JdHP}}
We mimick the proof of Theorem 4.3 and 4.4 of \cite{FGAR} (see also \cite{GaNaRomat}).  
The weak convergence topology on $\mathcal M_1 (\mathbb T)$ is equivalent to the topology of convergence of moments on $\bar{\mathbb D}^{\mathbb N_0}$.

The sequence $\mu^{(n)}$ is exponentially tight since we work on $\mathbb T$.
The mapping 
\begin{align*}
m : \mathcal M_1(\mathbb T) \rightarrow \bar{\mathbb D}^{\mathbb N_0}, \qquad  m(\mu) : = \left(m_k(\mu) := \int_{\mathbb T} z^k d\mu(z)\right)_{k \geq 1} 
\end{align*}
being a continuous injection, the LDP on $\mathcal M_1(\mathbb T)$ is then a consequence of the following LDP on the sequence of moments and of the inverse contraction principle (see \cite{demboz98} Theorem 4.2.4 and the subsequent Remark (a)).

\begin{prop} 
The sequence $(m( \mu^{(n)}))_n$ satisfies the LDP in $\bar{\mathbb D}^{\mathbb N_0}$ with speed $n$ and good rate function $\mathcal I_{\operatorname{m}}$ defined as follows. This function is finite in $(m_1, m_2, \dots)$ if and only if this is the moment sequence of a nontrivial measure $\mu \in \mathcal M_1(\mathbb T)$ with deformed Verblunsky coefficients $(\gamma_0, \gamma_1, \dots)\in {\mathbb D}^{\mathbb N_0}$ satisfying
\[\sum_{k=0}^\infty H_\d (\gamma_k) < \infty . \]
In that case
\begin{equation}
\mathcal I_{\operatorname{m}}(m_1,m_2,\dots) = \sum_{k=0}^\infty H_\d (\gamma_k)\,.
\end{equation}
\end{prop}
\proof
By Lemma \ref{LDPalphaHP} for any fixed $k$,  $(\gamma_0^{(n)}, \dots \gamma_{k-1}^{(n)})_{n \geq k}$ 
satisfies the LDP in $\bar{\mathbb D}^k$  with good rate function
\[\mathcal I^{(k)}(\gamma_0, \dots, \gamma_{k-1}) = 
 \sum_{j=0}^{k-1} H_\d (\gamma_j)\,.\]
 By contraction, 
this yields the LDP for  $(\alpha_0^{(n)}, \dots , \alpha_{k-1}^{(n)})_{n \geq k}$ in  $\bar{\mathbb D}^k$, and then for the finite sequence of moments
$(m_1(\mu^{(n)}), \dots, m_{k}(\mu^{(n)}))$. The rate for the latter LDP is 
\begin{align*}
\mathcal I_{\operatorname{m}}^{(k)}(m_1, \dots, m_{k}) = \mathcal I^{(k)}(\gamma_0, \dots, \gamma_{k-1}) , 
\end{align*}
where $\gamma_0, \dots, \gamma_{k-1}$ are the uniquely determined first deformed Verblunsky coefficients of any measure with first moments $m_1,\dots ,m_k$. In particular, $\mathcal I_{\operatorname{m}}^{(k)}(m_1, \dots, m_{k}) =\infty$ if no such measure exists, or if this measure has less than $k$ support points.

By the projective method of Dawson-G\"{a}rtner's theorem (Theorem 4.6.1 in \cite{demboz98}), the sequence $(m( \mu^{(n)}))_n$ satisfies the LDP in $\bar{\mathbb D}^{\mathbb N_0}$ with speed $n$ and good rate function
\begin{align*}
\mathcal{I}_{\operatorname{m}}(m_1,\dots) &= \sup_{k\geq 1} {\mathcal{I}}_{\operatorname{m}}^{(k)}(m_1,\dots m_{2k-1})  .
\end{align*}
This supremum is infinite if $(m_1, m_2,\dots )$ is not the moment sequence of a nontrivial probability measure on $\mathbb{T}$. 
Otherwise, there exists a unique sequence of deformed Verblunsky coefficients $(\gamma_0, \gamma_1, \dots)\in {\mathbb D}^{\mathbb N_0}$ corresponding to this measure and 
\begin{align*}
\mathcal{I}_{\operatorname{m}}(m_1,\dots) = \sup_{k\geq 1} \mathcal I^{(k)}(\gamma_0, \dots, \gamma_{k-1}) =  \sup_{k\geq 1} \sum_{j=0}^k H_\d (\gamma_k) = \sum_{j=0}^\infty H_\d (\gamma_j). 
\end{align*}
\hfill $\Box$

\subsection{Proofs of Section \ref{sumrule}}

\subsubsection{Proof of Corollary \ref{coreasy}}
The elementary decomposition 
\[(1- \g\cos \theta) = \g (1- \cos \theta) + (1 - \g)\]
and the definition of $\GW_{-\g}$ give
\[\mathcal K(\GW_{-\g}| \mu) =  H(\g) +\g \mathcal K(\GW_{-\mathbf 1}| \mu) + (1-\g) \mathcal K(\GW_{\mathbf 0}| \mu)- \g H(1). \]
Where for $|a|<1$, 
\begin{equation}
H(a) :=\int_0^{2\pi}   (1 - a \cos \theta) \log (1 - a \cos \theta) \ \frac{d\theta}{2\pi}
=  1 - \sqrt{1- a^2} + \log \frac{1 + \sqrt{1-a^2}}{2}\,.
\end{equation}

\begin{rem}
The minimum in formula (\ref{sumrulegwg}) is $0$. It is reached uniquely at $\mu=\GW_{-\g}$ corresponding to the  Verblunsky coefficients given in
(\ref{alphalimGW}). 
\end{rem}
\subsection{Proofs of Section \ref{smamat}}
\subsubsection{Proof of Proposition \ref{identif}}

First, we have $\ll \boldsymbol{\varepsilon}_1 , \mathcal R_1 \boldsymbol{\varepsilon}_1  \gg \, = \boldsymbol{c}_0(\mathcal G)$ and for $k\geq 2$,
\begin{eqnarray*}
\ll \boldsymbol{\varepsilon}_k , \mathcal R_k \boldsymbol{\varepsilon}_k  \gg \!\ 
= \!\ \ll  \varepsilon_k , \mathcal R_{k-1}^\dagger\mathcal R_{k-2}^\dagger \dots \mathcal R_1^\dagger \mathcal G \boldsymbol{\varepsilon}_k \gg \!\ = \!\ \ll \mathcal R_1 \dots \mathcal R_{k-1} \boldsymbol{\varepsilon}_k , \mathcal G\boldsymbol{\varepsilon}_k \gg\\
\!\ = \!\ \ll \mathcal G\hat\pi_{k-1}(\mathcal G)^\dagger\boldsymbol{\varepsilon}_k , \mathcal G\boldsymbol{\varepsilon}_k\gg\!\ = \!\ \ll \hat\pi_{k-1}(U)^\dagger\boldsymbol{\varepsilon}_k , \boldsymbol{\varepsilon}_k \gg \!\  =\!\ \ll  \boldsymbol{\varepsilon}_k ,  \hat\pi_{k-1}(\mathcal G)\boldsymbol{\varepsilon}_k \gg\!\  = \mathbf{c}_{k-1} (\mathcal G)\,. \end{eqnarray*}

To compute $\mathbf{c}_k (\mathcal G)$ we start from the definitions of $\mathcal G$ and $\Xi_N^p$, which yield 
\begin{eqnarray}
\label{xiroot}
 \Xi_N^p (\mathcal G^R) = \hat\Theta (u_0) \mathcal G^R(\abold_1, \abold_2, \dots)
\end{eqnarray}
where, if $u \in \mathbb U(p)$
\[\hat\Theta(u) = \begin{pmatrix} u & 0_{p, N-2p}\\
0_{N-2p, p}& {\mathrm I}_{N-2p}
\end{pmatrix}
\]
and 
\[u_0 = - \abold_0 + \rbold_0^R (\mathbf{1}- \abold_0^\dagger)^{-1}\rbold_0^L\,,\]
so that
\[\mathbf{c}_1 (\mathcal G^R) = [ \Xi_N^p (\mathcal G^R)]_p =u_0 \abold_1^\dagger\,. \]
More generally, looking for a recursion - thinking of (\ref{projxi}) - , we notice that
\begin{eqnarray}
\label{thetag}
\Xi_{N}^p \left(\hat\Theta(u) \mathcal G^R (\abold_0, \dots, \abold_{n-1})\right) =
\hat \Theta(v)\mathcal G^R(\abold_1, \dots, \abold_{n-1})
\end{eqnarray} where
\begin{eqnarray}v = v(u, \abold_0) = -\abold_0 + \rbold_0^R (\mathbf{1} - u\abold_0^\dagger)^{-1}u \rbold_0^L\,.
\end{eqnarray}
We need the following result.
\begin{lem}
\label{obs}
If $\abold \in \mathbb B_p$ and  $u \in \mathbb U(p)$ then
\begin{equation}
- \abold + \rbold^R (\mathbf{1} - u\abold^\dagger)^{-1}u\rbold^L = (\rbold^R)^{-1} (u -\abold) (\mathbf{1} - \abold^\dagger u)^{-1}\rbold^L = \boldsymbol{T}_{\abold}(u)\,.
\end{equation} 
\end{lem}
Let us assume that
\begin{equation}
\label{HR}\Xi_N^{jp} = \hat\Theta(u_j) \mathcal G^R(\abold_j, \dots, \abold_{n-1})\,,\end{equation}
where $u_j$ depends on $\boldsymbol{\alpha}_0, \dots, \boldsymbol{\alpha}_{j-1}$. 

Applying (\ref{projxi}),  (\ref{thetag}) and Lemma \ref{obs}  we get
\begin{eqnarray*}\Xi_N^{(j+1)p} (\mathcal G^R(\abold_0, \dots, \abold_{n-1})) &=& \Xi_{N-jp}^p \left( \hat\Theta(u_j) \mathcal G^R(\abold_j, \dots, \abold_{n-1})\right)\\
&=& \hat\Theta\left(\mathbf{T}_{\abold_j}(u_j)\right)\mathcal G^R (\abold_{j+1}, \dots, \abold_{n-1})\,,
\end{eqnarray*}
and the assumption (\ref{HR}) is satisfied at rank $j+1$ with 
\begin{equation}
\label{recu}u_{j+1} = \mathbf{T}_{\abold_j}(u_j)\,.\end{equation}
Passing to the upperleft block, we obtain easily, for every $j \leq n-1$
\[\mathbf{c}_{j} 
= u_{j}\abold_{j}^\dagger\,.\]
Now, using (\ref{magicT}), we see that $\boldsymbol{\beta}_j :=u_j^{-1}$ satisfies the recursion
\[\boldsymbol{\beta}_{j+1} = \mathbf{T}_{\abold_{j}^\dagger}(\boldsymbol{\beta}_j)\]
or, reversing
\[\boldsymbol{\beta}_j = \mathbf{T}_{-{\abold_j}^\dagger}(\boldsymbol{\beta}_{j+1})\]
which allows to conclude that $\boldsymbol{\beta}_j = \bbold_j$ and ends the proof of Proposition \ref{identif}.

\medskip

\paragraph{Proof of Lemma \ref{obs}:}
We have to prove 
\begin{align*} 
\rbold^R(- \abold + \rbold^R (\mathbf{1} - u\abold^\dagger)^{-1}u\rbold^L)(\rbold^L)^{-1} 
= (u -\abold) (\mathbf{1} - \abold^\dagger u)^{-1} .
\end{align*}
Since $\rbold^R \abold = \abold \rbold^L$ and $(\rbold^R)^2 = \mathbf{1} - \abold\abold^\dagger$, the left hand side simplifies to  
\begin{align*}
 - \abold + (\mathbf{1} - \abold\abold^\dagger) (\mathbf{1}- u\abold^\dagger)^{-1}u\,.
\end{align*}
Now, $(\mathbf{1} - u \abold^\dagger)u = u (\mathbf{1} -\abold^\dagger u)$ such that the last line is equal to 
\[  - \abold \left(\mathbf{1} + \abold^\dagger u (\mathbf{1}- \abold^\dagger u)^{-1}\right) + u (\mathbf{1}- \abold^\dagger u)^{-1}\,,\]
which is exactly $(u-\abold)(\mathbf{1} - \abold^\dagger u)^{-1}$.  \hfill $\Box$

\subsubsection{Proof of Theorem \ref{alphajlaw}} 
\label{sec:alphajlaw}
In  \cite{FGARop} it is proved that when $\mathbb U(N)$ is equipped with the Haar measure $\mathbb P^{(N)}$, the distribution of $(\abold_0, \cdots \abold_{n-1})$ is up to a normalization constant
\begin{equation}
\label{lawyA}
 \left(\otimes_{r= 0}^{n-2}( \det(\mathbf{1} - \abold_r \abold_r^\dagger)^{N- (r+2)p} \ 
d\abold_r) \right)\otimes
\ d\mathbb P^{(p)}(\abold_{n-1})\,,
\end{equation}
where, for $r=0, \dots, n-2$, $d\abold_r$ denotes the Lebesgue measure on $\mathbb B_p$, and 
$\mathbb P^{(p)}$ is, as usual, 
the Haar measure on $\mathbb U(p)$.

Since $\gbold_r$ is  $\abold_r^\dagger$ up to multiplication by a unitary matrix depending only on $(\abold_0, \cdots, \abold_{r-1})$, we deduce that, the pushforward of $\mathbb P^{(N)}$ 
 by  
 $(\gbold_0, \gbold_1, \cdots, \gbold_{n - 1})$ has again the distribution proportional to
\begin{equation}
\label{lawyH}
 \left(\otimes_{r= 0}^{n-2} (\det(\mathbf{1} - \gbold_r^\dagger \gbold_r)^{N- (r+2)p} \ 
d\gbold_r)\right)\otimes  d\mathbb P^{(p)}(\gbold_{n-1})\,,
\end{equation}

Now, by definition
\[\frac{d\mathbb H\mathbb P_\delta^{(N)}}{d \mathbb P^{(N)}}(U) = \hbox{const} \cdot \ \det(\mathrm{I}_N - U)^{\bar\delta} \det(\mathrm{I}_N - U^\dagger)^\delta\,.\]
It remains to apply (\ref{detfond}) and Proposition \ref{identif} to conclude that 
under $\mathbb H\mathbb P_\delta^{(N)}$, the variables $(\gbold_0, \gbold_1, \cdots, \gbold_{n - 1})$ are  independent and for $0 \leq r \leq n-2$ the density of $\gbold_r^p$ in $\mathbb B_p$  is proportional to
\[ \left(\det(\mathbf{1} - \gbold)\right)^{\bar \delta}\left(\det(\mathbf{1} - \gbold^\dagger)\right)^{ \delta} \det( \mathbf{1} - \gbold^\dagger \gbold)^{N- (r+2)p}\,.  \]
Further, the variable $\gbold_{n-1}$ has the distribution $\mathbb H \mathbb P_\delta^{(p)}$ on $\mathbb U(p)$. 
The value of the normalizing constant (\ref{cstner}) is then taken from formula (2.9) in \cite{Ner1}. \hfill $\Box$

\begin{rem}
Theorem 1.3 of Neretin \cite{Ner1} says that  if $\mathbb{U}(N)$ is equipped with the Haar measure, then  the distribution of $(\cbold_0 (U) , \dots, \cbold_{n-1}(U))$ is also (\ref{lawyH}). From (\ref{capital}) we deduce that, under $\mathbb H\mathbb P_\delta\sn$, 
$(\cbold_0 (U) , \dots, \cbold_{n-1}(U))$ and $(\gbold_0, \gbold_1, \cdots, \gbold_{n - 1})$ have the same distribution. 
The difference is that the second array depends only on the spectral measure, and the first one depends more deeply on $U$. In particular,  we do not know  the connection between these coefficients $\cbold(U)$ and $\boldsymbol{\alpha}$. 
\end{rem}

\subsubsection{Proof of Proposition \ref{lemmat}} \label{sec:proofoflemmat}
By independence, it suffices to prove the LDP only for one $\gbold_j^{(n)}$ with rate $H_{\d, p}$. Since the LDP is a standard consequence of the explicit density in \eqref{densalphajmat}, we only give a sketch of the proof. 
First, we get from the explicit expression of the constant in \eqref{cstner} 
\begin{align*}
\lim_{n\to \infty} \frac{1}{np} \log K_{n,j}^{(n\d)} = p H_\d (0) .
\end{align*}
Then, on the set $\{M \in \mathbb{C}^{p\times p}|\, MM^\dagger< \mathbf 1\}$ the rate function is finite and continuous. Indeed, if $\gbold \in \mathbb{C}^{p\times p}$ is a matrix with singular values smaller than 1, then $\mathbf 1-\gbold$ is non-singular. On the other hand, if $\gbold \in \mathbb{B}_p \setminus \{M \in \mathbb{C}^{p\times p}|\, MM^\dagger< \mathbf 1\}$, we have $H_{\d, p}(\gbold)=\infty$. This implies for any $\gbold\in \mathbb B_p$, denoting by $B_\varepsilon(\gbold)$ the open ball centered at $\gbold$ with radius $\varepsilon$ in the Frobenius norm, that
\begin{align*}
\lim_{\varepsilon\to 0} \limsup_{n\to \infty} \frac{1}{np} \log \mathbb H\mathbb P_{N\d}^{(N)} (\gbold_j^{(n)} \in B_\varepsilon(\gbold)) & = H_{\d, p}(\gbold), \\
\lim_{\varepsilon\to 0} \liminf_{n\to \infty} \frac{1}{np} \log \mathbb H\mathbb P_{N\d}^{(N)} (\gbold_j\sn \in B_\varepsilon(\gbold)) & = H_{\d, p}(\gbold) .
\end{align*}
From these limits, we get that $(\gbold_j\sn)$ satisfies the weak LDP with speed $N=np$ and good rate function $H_{\d, p}$. Necessarily, this sequence is exponentially tight, since it lives on the compact set $\mathbb B_p$, and the full LDP follows. \hfill $\Box$

\section*{Acknowledgments}
We warmly thank Barry Simon for his helpful remarks and suggestions about the gems, and Ofer Zeitouni for valuable conversations.
We also thank the anonymous reviewers for their careful reading of our manuscript and their many    
insightful comments and suggestions.  
%

\end{document}